\documentclass{amsart}

\usepackage{enumitem}
\usepackage[all]{xy}
\usepackage[margin=1.00in, marginparwidth=1.5cm, marginparsep=0.5cm]{geometry}
 \usepackage{booktabs}
 \usepackage{siunitx}
\usepackage[table]{xcolor} 
\usepackage{mathtools}
\usepackage{bbold}
\usepackage{placeins}
\usepackage{amssymb}
\usepackage{array}
\usepackage{xcolor}
\usepackage{mathrsfs}
\usepackage{physics}
\usepackage{amsthm}
\usepackage{amsmath,amsfonts,graphicx,mathdots,cite,wrapfig}
\usepackage[colorlinks,linkcolor=blue,citecolor=blue,pagebackref,hypertexnames=false, breaklinks]{hyperref}
\usepackage{cite}
\usepackage{mwe}
\usepackage{caption}
\usepackage{caption} 
\usepackage{subfigure} 
\usepackage{marginnote}

\newtheorem{theorem}{Theorem}[section]
\newtheorem{lemma}[theorem]{Lemma}
\newtheorem{proposition}[theorem]{Proposition}
\newtheorem{corollary}[theorem]{Corollary}
\newtheorem{assumption}[theorem]{Assumption}
\theoremstyle{definition}
\newtheorem{definition}[theorem]{Definition}

\theoremstyle{remark}
\newtheorem{remark}[theorem]{Remark}

\numberwithin{equation}{section}

\DeclareMathOperator{\jacobi}{\textbf{J}}
\DeclareMathOperator{\rr}{\mathbb{R}}
\DeclareMathOperator{\cc}{\mathbb{C}}
\DeclareMathOperator{\complex}{\mathbb{C}}
\DeclareMathOperator{\integers}{{\mathbb{Z}}}
\DeclareMathOperator{\nn}{\mathbb{N}}
\newcommand{\DiriAMO}{\mbox{\normalfont\Large\bfseries $\jacobi_{cs}^{D}$}}

\begin{document}

%%%%%%%%%%%%%%%%%%

%%%%%%%%%%%%%%%%%%

\title{Spectral decimation of piecewise centrosymmetric Jacobi operators on graphs}

%%%%%%%%%%%%%%%%%%

%    Information for first author
\author{Gamal Mograby}
\address{Department of Mathematics, University of Maryland, College Park, MD 20742, USA}
%    Current address
%\curraddr{Department of Mathematics and Statistics,
%Case Western Reserve University, Cleveland, Ohio %43403}
\email{gmograby@umd.edu, gamal.mograby@uconn.edu}
%    \thanks will become a 1st page footnote.
\thanks{The authors are very grateful to Robert Strichartz who, although no longer with us, continues to inspire research in harmonic and functional analysis, spectral theory, fractals and self-similarity. }
\dedicatory{In memory of Robert Strichartz}

%%    Information for second author
\author{Radhakrishnan Balu}
\address{Radhakrishnan Balu, Department of Mathematics \& Norbert Wiener Center for
	Harmonic Analysis and Applications, University of Maryland, College Park, MD 20742, USA}
\email{rbalu@umd.edu}
%%\thanks{Support information for the second author.}

%%    Information for second author
\author{Kasso A. Okoudjou}
\address{Kasso A. Okoudjou, Department of Mathematics, Tufts University, Medford, MA 02155, USA}
\email{kasso.okoudjou@tufts.edu}
%\thanks{Partially supported by the National Science Foundation under Grant No. DMS-1814253 and by ARO grant W911NF1910366}

%%    Information for second author
\author{Alexander Teplyaev}
\address{Alexander Teplyaev, Department of Mathematics, University of Connecticut, Storrs, CT 06269, USA}
\email{alexander.teplyaev@uconn.edu}
%\thanks{Partially supported by NSF DMS grant 1613025}

%    General info
\subjclass[2010]{47B36, 28A80, 81T17, 47A10, 05C20, 05C22, 81Q35}

\date{\today}

\keywords{Spectral Analysis, Jacobi operators, Graphs, Spectral Decimation, Renormalization group method}

\begin{abstract}
We study the spectral theory of a class of piecewise centrosymmetric Jacobi operators defined on an associated family of substitution graphs. Given a finite centrosymmetric matrix viewed as a weight matrix on a finite directed path graph and a probabilistic Laplacian viewed as a weight matrix on a locally finite strongly connected graph, we construct a new graph and a new operator by edge substitution. Our main result proves that the spectral theory of the piecewise centrosymmetric Jacobi operator can be explicitly related to the spectral theory of the probabilistic Laplacian using certain orthogonal polynomials. Our main tools involve the so-called spectral decimation, known from the analysis on fractals, and the classical Schur complement.  We include several examples of self-similar Jacobi matrices that fit into our framework.
\end{abstract}

\maketitle

\tableofcontents

\section{Introduction}
%{\color{blue} 
We study the spectral theory of a class of Jacobi operators, following B. Simon et al. \cite{BSimon2020PeriodicJacobi,CSZ21}, defined on an associated  family of  substitution graphs.  
Given a finite centrosymmetric matrix viewed as a weight matrix on a finite directed path graph and a probabilistic Laplacian viewed as a weight matrix on a locally finite strongly connected direct graph, we construct a new (possibly high dimensional, large or infinite) graph and a new  operator that acts as a weight matrix on this graph. In a sense that will be made precise, this  graph can be viewed as a result of a substitution procedure involving the initial graphs, and as such will be referred to as a substitution graph. Further, the new operator, which we call a piecewise centrosymmetric Jacobi operator acts  on the space of square summable functions defined on the vertices of the substitution graph. The goal of this paper is to study the spectral theory of these piecewise centrosymmetric Jacobi operators. The construction of the substitution graphs allows us to  show that these spectra are related to the spectra of the initial finite centrosymmetric matrices and  can be determined using the spectral decimation method popularized in mathematical physics and analysis on fractals beginning with work of Rammal and Toulouse \cite{RammalToulouse1983,Rammal1984}. 
This method is equivalent to the   classical Dirichlet to Neumann map and Schur complement, see  \cite{KZ19,BK20,BKbook,EKKST,A-RCRST,BajorinVibration3Ngasket2008,MalozemovTeplyaev2003} and references therein.
We make an explicit link between the classical spectral theory   and the spectral decimation method   through the analysis of   families of orthogonal polynomials.  We conclude the paper with several examples of one-dimensional graphs that fit into our framework. These self-similar graphs are obtained using a substitution method starting from a pair of initial graphs and this framework allows us to consider finite and infinite self-similar graphs defined in one or higher dimensions.  The Jacobi operators  are then defined so as to respect the adjacency relations in the resulting graph, and can be naturally viewed as  weight matrices on the graphs.  Thus  we say that a  matrix $\jacobi=\big( \jacobi(x,y) \big)$ is a Jacobi matrix if   $\jacobi$ reflects the adjacency relations of the graph $G$, where ${x,y \in V(G)}$ and $G=(V(G),E(G))$ is a graph.  In particular, for two different vertices $x,y \in V(G)$, we set $\jacobi(x,y)=0$, whenever  $x$ and $y$ are not adjacent, i.e.,  $(x,y) \notin E(G)$; see  \cite{BSimon2020PeriodicJacobi} for a   definition of Jacobi matrices on periodic trees.
%  In our setting, it becomes natural to view a Jacobi operator as a weight matrix of the underlying graph. %We will use this framework to establish a one-to-one correspondence between the Jacobi operator and the underlying graph.  
%}

We begin by previewing the substitution method used to define the graphs of interest to our work, and refer to Section~\ref{sec:SubOperatorSection}  for more details. We first introduce the two initial inputs of the method:  the \emph{building block graph  $G_{cs}$} and the \emph{model graph $G_p$}. 

\begin{enumerate}
\item[(a)]  {\bf The building block graph $G_{cs}$:}
 Let $\jacobi_{cs}$ be a finite tridiagonal matrix and assume that $\jacobi_{cs}$ is centrosymmetric, i.e.  $\jacobi_{cs} = R \jacobi_{cs} R$, where $R$ is an anti-diagonal identity matrix
\begin{align*}
R =  \left(\begin{matrix}
0 & 0 & \dots & 0 & 0 & 1\\
0 & 0 & \dots & 0 & 1 & 0\\
0 & 0 & \dots & 1 & 0 & 0\\
\vdots & \vdots &  & \vdots & \vdots & \vdots \\
0 & 1 & \dots & 0 & 0 & 0\\
1 & 0 & \dots & 0 & 0 & 0
\end{matrix}\right)
\end{align*}
%The subscript $cs$ stands for \textit{centrosymmetric}. 
We regard $\jacobi_{cs}$ as the weight matrix of a finite directed path graph $G_{cs}$, in the sense that the diagonal and off-diagonal entries in $\jacobi_{cs}$ respectively denote the weights assigned to the vertices and edges of $G_{cs}$. 
We note that centrosymmetric matrices have a rich eigenstructure  and appear in several applications  \cite{SolitonsGradedGraphs2021Mograby, PQST2020Mograby, SpectraPST2021Mograby, PropCentrosyMatrices2005Liu,InvEigenprobCentrosym2004BaiChan,spectralCentrosym2001TaoYasuda,Centrosym1998Andrew,Weaver1985Centrosymmetry,VeinDaleBookCentrosymmetry1999,CantoniButler1976centrosymmetric,centrosymmetric1962Collar}.
Figure~\ref{fig:InitialCopyOfJacobi} shows an example of a building block graph $G_{sc}$  associated with the centrosymmetric matrix $\jacobi_{cs}$. 
%\begin{figure}[htp]
%\begin{minipage}[b]{1.0\textwidth}
%\centering
%%\hspace*{-4cm} 
%\resizebox{!}{!}{\input{Figures/InitialCopyOfJacobi.tikz}}
%\end{minipage}%
%\caption{A finite directed weighted path graph $G_{cs}$.  The weights associated with the edges and vertices are the entries of a centrosymmetric Jacobi matrix  $\jacobi_{cs}$.}
%\label{fig:InitialCopyOfJacobi}
%\end{figure}	
\begin{figure}[htp]
\centering
%\hspace*{-4cm} 
 \includegraphics[width=1.\textwidth]{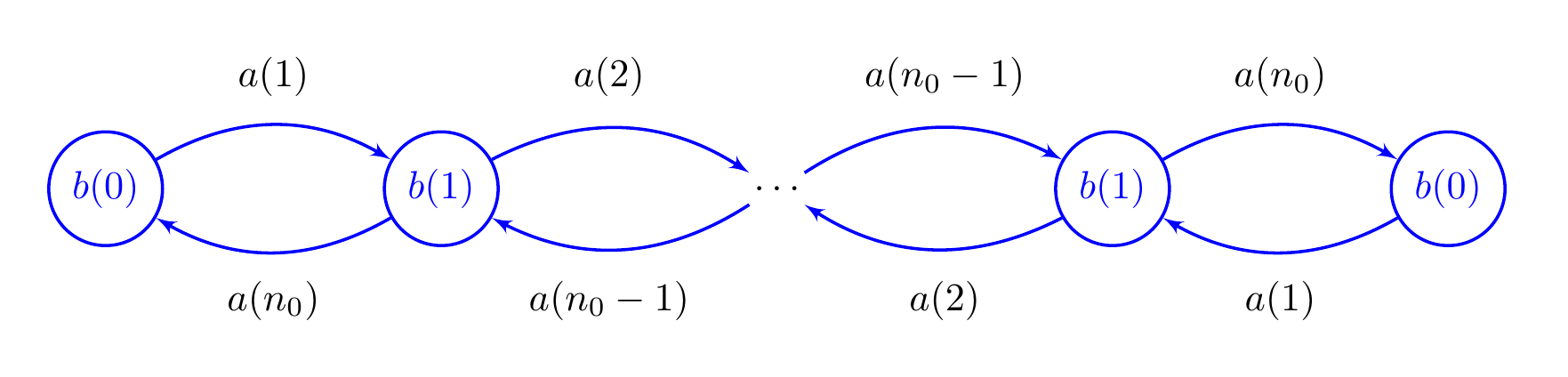}
\caption{A finite directed weighted path graph $G_{cs}$.  The weights associated with the edges and vertices are the entries of a centrosymmetric Jacobi matrix  $\jacobi_{cs}$.}
\label{fig:InitialCopyOfJacobi}
\end{figure}

\item[(b)]  {\bf The model graph  $G_p$:} Let $G_p$  be a locally finite connected directed graph. We consider a random walk on $G_p$
and therefore assume that the weight matrix of $G_p$ is given by a \textit{probabilistic graph Laplacian}
\begin{align*}
%\label{eq:probLapOnGp}
 \Delta_{\mathbf{p}} f (x)= f(x) \  - \sum_{y :(x,y) \in E(G_{p})} p(x,y) f(y),
\end{align*}
where  $\mathbf{p} = \{p(x,y)\}_{(x,y) \in E(G_{p})}$ is a given sequence of transition probabilities. Note that if we consider a symmetric random walk on $G_p$, then $ \Delta_{\mathbf{p}}$ becomes the standard probabilistic graph Laplacian. Throughout this work, we fix $G_p$ and assume that $ \Delta_{\mathbf{p}}$ is a  probabilistic graph Laplacian on $G_p$. 
%We refer to Section~\ref{sec:ProbabilisticgraphLaplacian} and Definition \ref{def:probLapOnGp} for details. 
%See precise definitions  in \ref{eq:CodForTransProb} and Definition \ref{def:probLapOnGp}. 

 %The subscript $p$ stands for \textit{probabilistic graph Laplacians} as
 %the graph $G_{p}$ will provide such Laplacians in the discussion below (see the more precise statements in (\ref{eq:CodForTransProb}) and Definition \ref{def:probLapOnGp} below).
\end{enumerate}
Using the building block graph $G_{cs}$ and the model graph $G_p$ we obtain the substitution graphs as follows. 

\begin{enumerate}
\item[(c)] {\bf The substitution graphs:}
Given a building block $G_{cs}$ and a model graph $G_p$, we construct a new graph $G=(V(G),E(G))$ by substituting copies of  $G_{cs}$ between adjacent vertices in $G_{p}$. Figure~\ref{fig:SubstitJacobiSG} shows an example of a substitution graph starting from the building block graph $G_{cs}$ of Figure~\ref{fig:InitialCopyOfJacobi}. For a choice of the probabilistic Laplacian $\Delta_{\mathbf{p}}$ of $G_p$, the weight matrix of the newly constructed graph $G$ is given by a \textit{substitution operator}
\begin{align*}
(\Delta_{p}, \jacobi_{cs})  \longmapsto F_{G_{p}}(\Delta_{\mathbf{p}}, \jacobi_{cs}),
\end{align*}
which will be formally introduce in Definition~\ref{def:Substitution Operator}.  The substitution  operator assigns to each pair $(\Delta_{\mathbf{p}}, \jacobi_{cs}) $ a Jacobi matrix on the graph $G$, which we denote by $\jacobi:=F_{G_{p}}(\Delta_{\mathbf{p}}, \jacobi_{cs}) $ and refer to as a \textit{piecewise centrosymmetric Jacobi operator} on $G$. 
That the substitution operator is well-defined is a consequence of  the symmetry assumption on $\jacobi_{cs}$.  In addition, as we will see later, this symmetry assumption is key in establishing many of our results. 
%The substitution  operator assigns each pair $(\Delta_{\mathbf{p}}, \jacobi_{cs}) $ to a Jacobi matrix on a graph $G$, which we denote by $\jacobi:=F_{G_{p}}(\Delta_{\mathbf{p}}, \jacobi_{cs}) $ and refer to as a \textit{piecewise centrosymmetric Jacobi operator} on $G$. One of our goals is to study the spectral theory of these operators. 

\begin{figure}[htp]
\centering
%\hspace*{-4cm} 
 \includegraphics[width=1.\textwidth]{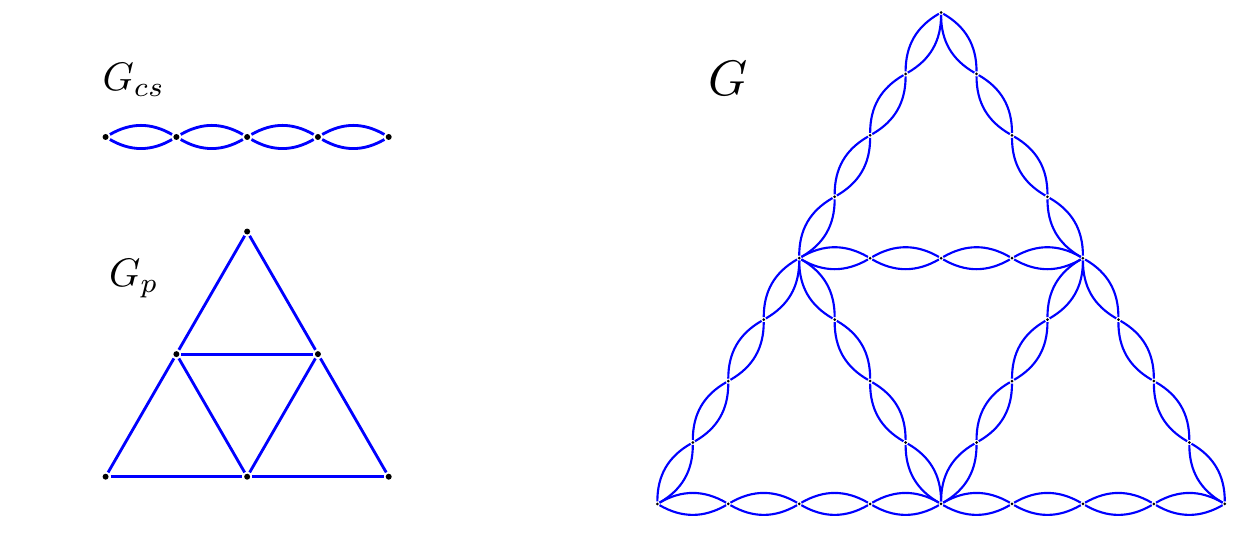}
\caption{ Let $G_{cs}$ be a finite directed weighted path graph associated with a centrosymmetric Jacobi matrix  $\jacobi_{cs}$ of Figure~\ref{fig:InitialCopyOfJacobi}. Let $G_p$ be a finite graph approximation of a Sierpinski lattice. We construct a new graph $G=(V(G),E(G))$ by substituting copies of the graph $G_{cs}$ between adjacent vertices in $G_{p}$.  The weight matrix of the graph $G$ will be given by $\jacobi:=F_{G_{p}}(\Delta_{\mathbf{p}}, \jacobi_{cs})$ and is an example of a piecewise centrosymmetric Jacobi operator on $G$.}
\label{fig:SubstitJacobiSG}
\end{figure}

%\begin{figure}[!htb]
%    \centering
%    \begin{minipage}{.5\textwidth}
%        \centering
%        %\hspace*{-2cm} 
%       \resizebox{4cm}{!}{\input{Figures/SGandJacobi.tikz}}
%    \end{minipage}%
%    \begin{minipage}{0.5\textwidth}
%        \centering
%        %\hspace*{-2cm} 
%        \resizebox{7cm}{!}{\input{Figures/SGlevel4.tikz}}
%    \end{minipage}
%    	\caption{ Let $G_{cs}$ be a finite directed weighted path graph associated with a centrosymmetric Jacobi matrix  $\jacobi_{cs}$ of Figure~\ref{fig:InitialCopyOfJacobi}. Let $G_p$ be a finite graph approximation of a Sierpinski lattice. We construct a new graph $G=(V(G),E(G))$ by substituting copies of the graph $G_{cs}$ between adjacent vertices in $G_{p}$.  The weight matrix of the graph $G$ will be given by $\jacobi:=F_{G_{p}}(\Delta_{\mathbf{p}}, \jacobi_{cs})$ and is an example of a piecewise centrosymmetric Jacobi operator on $G$.}
%	\label{fig:SubstitJacobiSG}
%\end{figure}

\end{enumerate}

The spectral theory of the piecewise centrosymmetric Jacobi operators will constitute the main contribution of this paper. To this end, in Section~\ref{sec:OrthogonalPolynomials} we begin the spectral analysis of the centrosymmetric matrix  $\jacobi_{cs}$ and connects it to the analysis of a related family  of polynomials and rational functions. 
This is summarized in our first main result, 
 %As noted above, the symmetry assumption on $\jacobi_{cs}$ guarantees that the substitution operator is well-defined, but it is also essential is establishing some of our main results. For example, we rely on this assumption to prove 
 Theorem \ref{thm:SchurToPoly} which allows us to find a relatively simple formula for the Schur complement of $\jacobi_{cs}$. Consequently,   the spectral analysis of the piecewise centrosymmetric Jacobi operators  $\jacobi=F_{G_{p}}(\Delta_{\mathbf{p}}, \jacobi_{cs}) $ can be investigated in the realm of the \textit{spectral decimation method}, a standard framework in the analysis on fractals. 
Developed in the 1980s  \cite{BellissardRenormalizationGroup,Rammal1984, RammalToulouse1983, Alexander1984,KadanoffAlexander1983}, the spectral decimation method attracted considerable attention over the last decades \cite{Shirai2000,BobTransfOfspectra2010,Strichartz2001MethodOfAverages,ShimaFukushima1992,ShimaSierpinski1993,ShimapreSierpinski1991,Teplyaev1998,BobBook2006,BajorinVibrationSpectra2008,KronTeufl2004}. 
For an overview of some modern approach to the spectral decimation method, we refer to 
\cite{KigamiBook2001,StrichartzTeplyaev2012,BobFractafolds2003,MalozemovTeplyaev2003,BobBook2006}. 
 
 The framework for the spectral decimation method developed in \cite{MalozemovTeplyaev2003}  is the most suited for the piecewise centrosymmetric Jacobi operators. Using it along with  a set of assumptions (Assumption~\ref{assump:HilbertspacesAndIsometry}) we obtain one of our other main results, Theorem~\ref{coro:RenoEquAndSpectra}, which gives a description of the spectrum of $\jacobi=F_{G_{p}}(\Delta_{\mathbf{p}}, \jacobi_{cs}) $.  The essence of this result is that it relates the spectra $\sigma \big( \jacobi  \big)$, $ \sigma(\Delta_{\mathbf{p}})$, and  an  \textit{exceptional set}  $\mathscr{E}_{\jacobi_{cs}}$ that appears naturally in the spectral decimation method. The link between these sets is obtained using a polynomial, \emph{the spectral decimation function} $R_{\jacobi_{cs}}$ that we are able to explicitly calculated. Furthermore, we prove that the resolvent operators of $\jacobi$ and $\Delta_{\mathbf{p}}$ satisfy the renormalization group identity given in \eqref{eq:RenoEquation}. This result is related to \cite[Theorem 2.2]{BGM-1988} and Bellissard's work on quasicrystals \cite[Theorem 1]{BellissardRenormalizationGroup}, although we do not rely on them. However, we note  that the connections and differences bewteen our framework and Bellissard's are elaborated in \cite[Section 6]{BaluMogOkoTep2021spectralAMO}.

The rest of the paper is organized as follows. In Section~\ref{sec:SubOperatorSection}, we define the substitution graph and operator associated with a fixed graph $G_p$ and a building block graph $G_{cs}$. We show that each element in the range of the latter can be identified with a piecewise centrosymmetric Jacobi operator. Proposition~\ref{prop:PropertiesOfFstar}  summarizes some of the relevant properties of the substitution operator.  Section~\ref{sec:OrthogonalPolynomials} is divided into two parts. The first part deals with the spectral analysis of a centrosymmetric matrix $\jacobi_{cs}$, to which is associated two families of polynomials. Using these polynomials, we  compute  in Theorem~\ref{thm:SchurToPoly} the Schur complement of $\jacobi_{cs}$, and identify  the polynomial $R_{\jacobi_{cs}}(z)$ as  a spectral decimation function.  Subsequently, the second part of Section~\ref{sec:OrthogonalPolynomials} is devoted to an analysis of $R_{\jacobi_{cs}}(z)$. In particular, we relate the  preimages of the points $\{0,2\}$ under $R_{\jacobi{cs}}$ to the exceptional set $\mathscr{E}_{\jacobi_{cs}}$, determine the critical points of $R_{\jacobi_{cs}}(z)$ thereby establishing  the existence of the branches of the inverse $R^{-1}_{\jacobi_{cs}}$ in the domain $[0,2]$. 
%These properties are crucial considering the fact that $\sigma(\Delta_{\mathbf{p}}) \subset [0,2]$, see Remark \ref{Rem:specSubsetOfInterval}. 
These results are summarized in Theorem \ref{thm:propertiesOfSpecDeci}.  Section \ref{sec:DeciSection} contains the main results on the spectra of piecewise centrosymmetric matrices. More specifically, in Proposition~\ref{prop:SpectralSimilarity}  we prove the spectral similarity between  $\jacobi$ and  $\Delta_{\mathbf{p}}$ and determine the  exceptional set $\mathscr{E}_{\jacobi_{cs}}$. We then prove our main result, Theorem \ref{coro:RenoEquAndSpectra} which gives a complete description of the spectrum of the piecewise centrosymmetric Jacobi operator $\jacobi$.   Finally, in Section \ref{sec:OnedimPathGraph}, we illustrate our results by focussing on one-dimensional path graphs which corresponds to the classical Jacobi matrices.

Our work is part of a long term study of  mathematical physics on fractals and graphs,  more specifically, quantum Hall systems with AMO and  their topological quantum phases \cite{ADRS21,KassoStrichartz2005,KassoStrichartz2007,KassoSaloffCosteTeplyaev2008,Akkermans2009ComplexDim,Akkermans2012ThermoPhoton,Akkermans2012SpatialLogPeriodic,Akkermans2020AC_circuits,Akkermans2013quantumFieldsFractals,Dunne2012HeatKernels,
AlonsoRuiz2016Hanoi,HinzMeinert2020}, in which novel features of physical systems can be associated with the unusual spectral and geometric properties of fractals and graphs compared to smooth manifolds.

 \section{Weighted substitution with Centrosymmetric  Jacobi matrices} \label{sec:SubOperatorSection}
This section introduces  the main concepts of the paper, namely the definition of the substitution graph and operator.

\subsection{Centrosymmetric Jacobi Matrices}\label{sec:centrosymmetryJacobis}

Suppose $G_{cs}=(V(G_{cs}),E(G_{cs}))$ is a finite directed weighted path graph with vertices $V(G_{cs})$ and edges $E(G_{cs})$.  We assume that $G_{cs}$ consists of $n_0+1$ vertices for $n_0 \geq 1$, and denote the vertices by $V(G_{cs})=\{0,1,\dots, n_0\}$. The set of edges is then given by $E(G_{cs}) = \big\{(i,i+1) \ \rvert \ i \in  \{0, \dots, n_0-1\} \big\} \cup \big\{(i+1,i) \ \rvert \ i \in  \{0, \dots, n_0-1\} \big\}$. We describe the weights assigned to the vertices and edges by the following $(n_0+1) \times (n_0+1)$ centrosymmetric Jacobi matrix
\begin{align}
\label{eq:BasicJacobi}
\jacobi_{cs} :=
\begin{pmatrix}
b(0) & a(1) & 0 &   \dots & 0   \\
a(n_0) & b(1) & a(2) &   \dots & 0  \\
0 & a(n_0-1) & b(2) &    \ddots &  \vdots  \\
\vdots & \vdots & \ddots & \ddots &   a(n_0)  \\
0 & 0 & 0 & a(1) & b(0)  
\end{pmatrix}
\end{align}
The boundary and interior vertices of $G_{cs}$ are given by $\partial G_{cs} := \{0, n_0\}$ and $V(G_{cs}) \backslash \partial G_{cs} = \{1, \dots ,n_0-1\}$, respectively. In the sequel, the class of centrosymmetric matrices will be denoted by 
\begin{align}
\label{eq:centroSymMatrices}
\mathcal{CSJ} := \{ \ \jacobi_{cs}  \rvert    \jacobi_{cs}  \text{ is centrosymmetric of the form } (\ref{eq:BasicJacobi}), \text{ where } n_0 \geq 1 \ \}.
\end{align}

\subsection{Probabilistic graph Laplacians}\label{sec:ProbabilisticgraphLaplacian}

Suppose $G_{p} = (V(G_{p}),E(G_{p}))$ is a locally finite (strongly) connected directed graph. Let $\{p(x,y)\}_{(x,y) \in E(G_{p})}$ be a sequence of weights assigned to the directed edges. The edge $(x,y)$ points from the vertex $x$ to $y$ and we regard $p(x,y)$ as a transition probability from $x$ to $y$. We impose the following conditions
\begin{align}
\label{eq:CodForTransProb}
\begin{cases}
\  (x,y) \in E(G_{p}), \ \Leftrightarrow  \ 0<p(x,y) \leq1 \\
\  (x,y) \notin E(G_{p}), \ \Leftrightarrow  \ p(x,y) =0 \\
\ \sum_{y :(x,y) \in E(G_{p})} p(x,y) = 1, \ \forall \ x \in V(G_{p}).
\end{cases}
\end{align}
\begin{definition}
\label{def:probLapOnGp}
Let $\mathbf{p} = \{p(x,y)\}_{(x,y) \in E(G_{p})}$ be a given sequence of transition probabilities, i.e. $\mathbf{p}$ satisfies the conditions in (\ref{eq:CodForTransProb}). The probabilistic graph Laplacian on the graph $G_{p}$ associated to $\mathbf{p}$ is defined by
\begin{align}
\label{eq:probLapOnGp}
 \Delta_{\mathbf{p}} f (x)= f(x) \  - \sum_{y :(x,y) \in E(G_{p})} p(x,y) f(y),
\end{align}
where $f \in \ell (G_{p}) := \{f:V(G_{p}) \to \complex \}$ and $x \in V(G_{p})$.
\end{definition}
The collection of all  probabilistic Laplacians on  a fixed graph $G_{p}$ is denoted by 
\begin{align}
\label{eq:setOfProbLaplacians}
\mathcal{L}_{G_{p}}:= \big\{\Delta_{\mathbf{p}} \ \big\rvert \  \Delta_{\mathbf{p}} \text{ defined in } (\ref{eq:probLapOnGp}), \mathbf{p}=\{p(x,y)\}_{(x,y) \in E(G_{p})} \text{ satisfies } (\ref{eq:CodForTransProb}) \big\}
\end{align}
\begin{remark}
The third equation in the conditions (\ref{eq:CodForTransProb}) arises naturally in different contexts and applications \cite{DurrettBookProbability2019, GrigoryanIntroGraphs2018,KellerLenzBook2021,KellyBookReversibility2011}. It is related to the consistency condition in \cite[equation (1.4)]{BobTransfOfspectra2010}, where  Strichartz utilizes the electrical network interpretation of the weights. This condition is also found in \cite[Assumption 2.8]{SolitonsGradedGraphs2021Mograby}, in the framework of Toda-lattices on fractal-type graphs as one of the sufficient conditions leading to the existence of static soliton solutions on such graph. This condition also makes it a Markov chain that later helps us to apply the Kolomogorov’s condition to establish operators of interest as self-adjoint.
\end{remark}
Throughout the paper, we make the following assumption on $G_{p}$:
\begin{assumption}
\label{ass:EdgesInBothDirections}
We assume that if $(x,y) \in E(G_{p})$, then, so is $(y,x) \in E(G_{p})$ and we refer to $x,y \in V(G_{p})$ as adjacent vertices. Note that in general, the transition probabilities $p(x,y)$ and $p(y,x)$ are not equal.
\end{assumption}

\subsection{Weighted  edge substitution operator}\label{sec:WeightedSubstitutionOperator}

We construct a new graph $G=(V(G),E(G))$ by starting from $G_{p}$, then substituting copies of a graph $G_{cs}$ between adjacent vertices in $G_{p}$. 
%Both graphs are defined as in the previous section 
We refer to \cite{BobTransfOfspectra2010} for similar constructions. Now recall that the vertices and boundary vertices of $G_{cs}$ are given by $V(G_{cs})  = \{0, 1, \dots ,n_0\}$ and $\partial G_{cs} = \{0, n_0\}$, respectively.  Let $(x,y) \in E(G_{p})$ (hence $(y,x) \in E(G_{p})$) for some $x,y \in V(G_{p})$ . We replace the edges $(x,y)$ and $(y,x)$ in the graph $G_{p}$ with a copy of $G_{cs}$, by identifying $x$ with $0$ and $y$ with $n_0$. Equivalently, we can identify $x$ with $n_0$ and $y$ with $0$, due to the centrosymmetry assumption on $G_{cs}$ and in either case, the resulting graph $G$ is the same. We distinguish two types of vertices in the graph $G$. The vertices which we obtain from $G_{p}$ and therefore being just $V(G_{p})$ and the vertices that are obtained from $G_{cs}$ through the substitution procedure. The edges of the graph $G$ are exactly the edges in each copy of $G_{cs}$. We construct Jacobi matrix on a graph $G$, and view it as  the graph weight matrix. We refer to \cite{BSimon2020PeriodicJacobi} for a similar definition.
\begin{definition}
A Jacobi matrix on a graph  $G=(V(G),E(G))$ is a matrix $\jacobi=\big( \jacobi(x,y) \big)_{x,y \in V(G)}$ indexed by the vertices $V(G)$, such that $\jacobi(x,y) = 0$ whenever $ (x,y) \notin E(G)$ and $x \neq y$. We denote the set of all Jacobi matrices on graphs by 
\begin{align}
\label{eq:JacobiMatrices}
\mathcal{J}:= \{ \ \jacobi  \rvert   \jacobi  \text{ is a Jacobi matrix on a graph} \}.
\end{align}
\end{definition}
\begin{definition}[Substitution Operator]
\label{def:Substitution Operator}
Fix a graph $G_{p}$ as defined in Section \ref{sec:ProbabilisticgraphLaplacian}. A substitution operator associated with $G_{p}$ is the mapping $F_{G_{p}}$ given by 
\begin{align*}
F_{G_{p}}: \mathcal{L}_{G_{p}}   \times \mathcal{CSJ} \to \mathcal{J} \\
(\Delta_{\mathbf{p}}, \jacobi_{cs})  \longmapsto F_{G_{p}}(\Delta_{\mathbf{p}}, \jacobi_{cs})
\end{align*}
where  $F_{G_{p}}(\Delta_{\mathbf{p}}, \jacobi_{cs})$ is defined as follows.  Let $G_{cs}$ be the graph associated to the centrosymmetric matrix $\jacobi_{cs}$ and $G=(V(G),E(G))$ is the resulting substitution graph. The substitution procedure naturally induces a \textit{covering map}, which we denote by $\varphi: V(G) \to V(G_{cs})$ (each vertex in $G$ is in a copy of $G_{cs}$ and hence corresponds naturally to a vertex in  $G_{cs}$):
\begin{enumerate}
\item $J:=F_{G_{p}}(\Delta_{\mathbf{p}}, \jacobi_c) $ is a Jacobi matrix on the graph $G$. The diagonal entries are given by
\begin{align*}
 \jacobi(x,x)   :=   \jacobi_{cs}(\varphi(x),\varphi(x)), \quad  \quad x \in V(G).
\end{align*}
\item Each $(x,y) \in E(G)$ is associated with an edge $(\varphi(x),\varphi(y))$ in a copy of  $G_{cs}$ and suppose that this copy replaced the edges between the vertices, say $u,v \in V(G_{p})$. Then we set
\begin{align}
\label{eq:gluingFormula}
\begin{cases}
\ \jacobi(x, y) := p(u,v)  \jacobi_{cs}(0,1), \quad \text{ if } x=u,\\
\ \jacobi(x, y) := \jacobi_{cs}(\varphi(x),\varphi(y)), \quad \text{ if } x\notin V(G_p).
\end{cases}
\end{align}
\end{enumerate}
\end{definition}
\begin{remark}
\label{rem:gluingProbLap}
In the particular case $n_0=1$, we have $V(G_{p}) = V(G)$. Moreover, if we choose $\Delta_{0}$ for the Jacobi matrix $\jacobi_{cs}$, where 
\begin{align}
\label{eq:TrivialProbLap}
 \Delta_{0} :=
\begin{pmatrix}
1 & -1   \\
-1 & 1
\end{pmatrix}.
\end{align}
then the Substitution formula (\ref{eq:gluingFormula}) gives $\jacobi = \Delta_{\mathbf{p}}$. This implies that the graphs $G$ and $G_{p}$ are identical.
\end{remark}
\begin{definition}
Let $\jacobi \in \mathcal{J}$. If there exist a graph $G_{p}$, and a probabilistic Laplacian $\Delta_{\mathbf{p}} \in \mathcal{L}_{G_{p}}$ and $ \jacobi_c \in \mathcal{CSJ}$, such that $\jacobi=F_{G_{p}}(\Delta_{\mathbf{p}}, \jacobi_c)$, then we call $\jacobi $ a piecewise centrosymmetric Jacobi operator.
\end{definition}
We end this section by summarizing the key properties of the substitution operator that will be needed in the sequel.

\begin{proposition}  The following statements hold.
\label{prop:PropertiesOfFstar}
\begin{enumerate}
	\item $F_{G_{p}}(\Delta_{\mathbf{p}}, \Delta_{0}) = \Delta_{\mathbf{p}}$, (where $\Delta_{0}$ is defined in (\ref{eq:TrivialProbLap})).
	\item $F_{G_{p}}(\Delta_{0}, \jacobi_{cs}) =  \jacobi_{cs}$.
	\item $F_{G_{p}}(\Delta_{\mathbf{p}}, I)=I $, where the identity matrices $I$ might be of different size.
	\item Let $\jacobi_{cs}^{(1)}$ and $\jacobi_{cs}^{(2)}$ be elements of $\mathcal{CSJ}$ and of the same size. Then for $\lambda, \mu \in \complex$, we have 
	\begin{equation*}
	F_{G_{p}}(\Delta_{\mathbf{p}}, \lambda \jacobi_{cs}^{(1)} + \mu \jacobi_{cs}^{(2)})  = \lambda F_{G_{p}}(\Delta_{\mathbf{p}}, \jacobi_{cs}^{(1)})  + \mu F_{G_{p}}(\Delta_{\mathbf{p}}, \jacobi_{cs}^{(2)}) .
	\end{equation*}
\end{enumerate}
\end{proposition}
\begin{proof}
The statements follow directly from the formula (\ref{eq:gluingFormula}).
\end{proof}

\section{Orthogonal polynomials and centrosymmetric  Jacobi matrices }\label{sec:OrthogonalPolynomials}
In this section we elaborate on the connection between the centrosymmetric Jacobi matrices and a family of orthogonal polynomials. 
\subsection{Schur complement of centrosymmetric Jacobi matrices}\label{sec:SchurcentrosymmetryJacobis}

The connection between orthogonal polynomials and the spectral theory proved fruitful during the last few decades and generated considerable interest among various communities of researchers. Our reference list is by no means exhaustive, and we refer the reader for instance to Simon’s treatises and the references therein \cite{SpectralTheoryOrthoPoly2014Simon,BrownEvansLittlejohn1992,EverittLittlejohn1991OrthoPoly,Uvarov1970, Berezansky2009,Simon2011SzegoTheorem,Simon2005OrthoPolyUnitCircleI}. The objective of this section is to reinterpret concepts from the spectral decimation technique in terms of orthogonal polynomials. Because we are dealing with Jacobi matrices,  orthogonal polynomials naturally arise through the associated three-term recurrence relations. For us, two cases are relevant, namely the three-term recurrence relations corresponding to the centrosymmetric  Jacobi matrices $\jacobi_{cs}$ and the restriction of $\jacobi_{cs}$ to the interior vertices of $G_{cs}$, which we denote by $\jacobi_{cs}^{D}$. Note that when dealing with $\jacobi_{cs}^{D}$ we assume $n_0 \geq 2$ as $\jacobi_{cs}^{D}$ is the interior block matrix of $\jacobi_{cs}$,

\begin{align}
\label{eq:DiriJacobi}
\jacobi_{cs}= 
\left( 
\begin{array}{c | c | c} 
b(0) & \begin{array}{c c c} 
     a(1) & 0 & \dots 
  \end{array} & 0 \\
 \hline 
\vdots &  \begin{array}{c c c} 
      &  &    \\ 
      & \DiriAMO & \\ 
      &  &  
  \end{array} & \vdots \\ 
  \hline 
0 &  \begin{array}{c c c} 
     \dots & 0 & a(1)
  \end{array} &  b(0)
 \end{array} 
\right).
 \end{align}
 \\
% \marginpar{(?) Kasso, the def of $\ell(I)$ is moved to the beginning of section}
 We regard $\jacobi_{cs}$ as a matrix acting on $ \ell \big( V(G_{cs})\big)=\ell \big( \{0,1,\dots, n_0\} \big)$, where $\ell (I)$ denotes the linear space of $\complex$-valued functions $\ell (I) := \{f:I \to \complex \}$ for a set $I$. In particular, $\jacobi_{cs}^{D}$ acts on $\ell \big( \{1,\dots, n_0-1\} \big)$. Throughout this section, we suppose that $\jacobi_{cs} \in \mathcal{CSJ}$ is fixed for some $n_0 \geq 2$ and $a(k)\neq 0 $ for each $k=0, \hdots, n_0+1$. We introduce a sequence of monic polynomials $P_0(z), P_1(z), \dots, P_{n_0+1}(z)$ corresponding to $\jacobi_{cs}$ as follows. We initialize $P_0(z):=1$ and $P_1(z):=z-b(0)$. For $k \in \{2,\dots, n_0+1\}$, we define $P_k(z) $ to be the determinant of the leading principal $k \times k$ submatrix of $zI-\jacobi_{cs}$, i.e.
 
{\fontsize{9}{10}
\begin{align}
P_k(z) 
:= \textbf{det}
\begin{pmatrix}
z-b(0) & -a(1) & 0 &   \dots & 0   \\
-a(n_0) & z-b(1) & -a(2) &   \dots & 0  \\
0 & -a(n_0-1) & z-b(2) &    \ddots &  \vdots  \\
\vdots & \vdots & \ddots & \ddots &   -a(k-1)  \\
0 & 0 & 0 & -a(n_0+2-k) & z-b(k-1) 
\end{pmatrix}.
\end{align} }%
\\
Note that $P_{n_0+1}(z) = \textbf{det} (zI-\jacobi_{cs})$ is the characteristic polynomial of $\jacobi_{cs}$   whose   eigenvalues  are the zeros of $P_{n_0+1}(z)$. It can be easily shown that the sequence of polynomials $P_0(z), P_1(z), \dots, P_{n_0+1}(z)$ satisfies the following recurrence relations
{\fontsize{8.5}{10}
\begin{align}
\label{eq:recurrenceRelations}
\begin{cases}
    P_0(z)=1, \quad P_1(z) = z-b(0)\\
    P_k(z) = \big(z-b(k-1)\big) P_{k-1}(z) - a(k-1)a(n_0+2-k)P_{k-2}(z) , \quad k \in \{2,\dots, n_0+1\}.
  \end{cases}
\end{align} }%
Similarly, for the Jacobi matrix $\jacobi^D_{cs}$, we initialize $P^D_0(z):=1$ and $P^D_1(z):=z-b(1)$. For $k \in \{2,\dots, n_0-1\}$, we define $P^D_k(z) $ to be the determinant of the leading principal $k \times k$ submatrix of $zI-\jacobi^D_{cs}$, i.e.
{\fontsize{9}{10}
\begin{align}
P^D_k(z) 
:= \textbf{det}
\begin{pmatrix}
z-b(1) & -a(2) & 0 &   \dots & 0   \\
-a(n_0-1) & z-b(2)  & -a(3) &   \dots & 0  \\
0 & -a(n_0-2) & z-b(3) &    \ddots &  \vdots  \\
\vdots & \vdots & \ddots & \ddots &   -a(k)  \\
0 & 0 & 0 & -a(n_0+1-k) & z-b(k) 
\end{pmatrix},
\end{align} }%
\\
where in this case  $P^D_{n_0-1}(z) = \textbf{det} (zI-\jacobi^D_{cs})$ is the characteristic polynomial of $\jacobi^D_{cs}$.
The sequence $P^D_0(z)$, $P^D_1(z)$, $\dots$, $P^D_{n_0-1}(z)$  satisfies the following recurrence relations
{\fontsize{9}{10}
\begin{align}
\label{eq:recurrenceRelationsDiri}
\begin{cases}
    P^D_0(z)=1, \quad P^D_1(z) = z-b(1)\\
    P^D_k(z) = \big(z-b(k)\big) P^D_{k-1}(z) - a(k)a(n_0+1-k)P^D_{k-2}(z) , \quad k \in \{2,\dots, n_0-1\}.
  \end{cases}
\end{align} }%
We first prove some auxiliary Lemmas.
\begin{lemma} \label{lem: helpfulIdentities1}
The following identities hold:
\begin{enumerate}
\item[(1)] $P^D_{1}(z)  P_{1}(z) - P^D_{0}(z) P_{2}(z) = a(1)a(n_0)$.
\item[(2)] For $n_0 \geq 2$ and $k \in \{2, \dots, n_0\}$, we have {\fontsize{9}{10}
\begin{equation*}
P^D_{k}(z)  P_{k}(z) - P^D_{k-1}(z) P_{k+1}(z)= a(n_0+1-k)a(k) \Big( P^D_{k-1}(z)  P_{k-1}(z) - P^D_{k-2}(z) P_{k}(z) \Big)
\end{equation*} }%
\item[(3)] {\fontsize{9}{10} $ P^D_{n_0-1}(z)  P_{n_0-1}(z) - P^D_{n_0-2}(z) P_{n_0}(z)= a(n_0)a(1) \Big( \prod^{n_0-1}_{i=2}a(i)  \Big)^2 $. }%
\end{enumerate}
\end{lemma}
\begin{proof} 
 (1) We use the recurrence relations (\ref{eq:recurrenceRelations}), (\ref{eq:recurrenceRelationsDiri})  and obtain {\fontsize{9}{10}
\begin{align*}
P^D_{1}(z)  P_{1}(z) - P^D_{0}(z) P_{2}(z)&= \big( z-b(1) \big)  P_{1}(z) - P_{2}(z) \\
%&=  a(1)a(n_0)P_{0}(z)  \\
&=  a(1)a(n_0)
\end{align*} }%
\begin{enumerate} 
\item[(2)] Similarly, repeated application of the recurrence relations (\ref{eq:recurrenceRelations}), (\ref{eq:recurrenceRelationsDiri}) gives, {\fontsize{9}{10}
\begin{align*}
P^D_{k}(z)  P_{k}(z) - P^D_{k-1}(z) P_{k+1}(z)= \\
=P^D_{k}(z)  P_{k}(z) - P^D_{k-1}(z) \Big(\big(z-b(k)\big) P_{k}(z) - a(k)a(n_0+1-k)P_{k-1}(z) \Big) \\
= \Big( P^D_{k}(z)  - P^D_{k-1}(z) \big(z-b(k)\big) \Big) P_{k}(z) + a(k)a(n_0+1-k)P_{k-1}(z)  P^D_{k-1}(z) \\
=  - a(k)a(n_0+1-k)P^D_{k-2}(z)P_{k}(z) + a(k)a(n_0+1-k)P_{k-1}(z)  P^D_{k-1}(z).
\end{align*} }%
\item[(3)] Iterate the arguments in parts (1) and (2).  
\end{enumerate}
\end{proof}
\begin{lemma}
\label{Lem:PhiNenner}
Let $n_0 \geq 2$. Then $$\big(P_{n_0}(z)\big)^2-\big( \prod^{n_0}_{i=1}a(i)  \big)^2 = \mathbf{det} \big(zI - \jacobi_{cs} )  \mathbf{det}(zI-\jacobi^D_{cs}).$$
\end{lemma}
\begin{proof}
We first compute $P_{n_0}(z)$ as a determinant and expand by the first row. Using the centrosymmetric assumption of $\jacobi^D_{cs}$, we see that the cofactor of $a(1)$ is $- a(n_0)P^D_{n_0-2}(z)$, leading to $P_{n_0}(z) = \big( z-b(0) \big) \mathbf{det} \big(zI-\jacobi^D_{cs} ) - a(1)a(n_0)P^D_{n_0-2}(z)$.
Hence
\begin{align}
\label{eq:identityComputDet1}
 \big( z- b(0) \big) P^D_{n_0-1}(z)   =& P_{n_0}(z) + a(1)a(n_0)P^D_{n_0-2}(z).
\end{align}
Repeated application of the recurrence relations (\ref{eq:recurrenceRelations}) and (\ref{eq:recurrenceRelationsDiri}) gives,{\fontsize{9}{10}
\begin{align*}
 \text{det} \big(zI - \jacobi_{cs} )  \text{det}(zI - \jacobi^D_{cs})=P_{n_0+1}(z) P^D_{n_0-1}(z) =\\
 = \big(z-b(0)\big) P_{n_0}(z) P^D_{n_0-1}(z)  - a(n_0)a(1) P_{n_0-1}(z) P^D_{n_0-1}(z) \\
 = \Big(  P_{n_0}(z) + a(1)a(n_0)P^D_{n_0-2}(z)   \Big) P_{n_0}(z) - a(n_0)a(1) P_{n_0-1}(z) P^D_{n_0-1}(z) \\
 %&= \big( P_{n_0}(z) \big)^2 + a(1)a(n_0)P^D_{n_0-2}(z) P_{n_0}(z) - a(n_0)a(1) P_{n_0-1}(z) P^D_{n_0-1}(z) \\
  = \big( P_{n_0}(z) \big)^2 - a(n_0)a(1) \Big( P^D_{n_0-1}(z) P_{n_0-1}(z) - P^D_{n_0-2}(z) P_{n_0}(z) \Big),
\end{align*}}%
where the third equality holds by (\ref{eq:identityComputDet1}). The statement follows by Lemma \ref{lem: helpfulIdentities1} (3).
\end{proof}

To state the main result of this section, we first establish some Schur complement computations for $\jacobi_{cs}$. These computations will be relevant in Section \ref{sec:DeciSection} and used in the context of the spectral decimation method. We decompose $\jacobi_{cs}$ on $ \ell \big( \partial G_{cs} \big) \oplus \ell \big( V(G_{cs}) \backslash \partial G_{cs} \big)=\ell \big( \{0,n_0\} \big) \oplus \ell \big( \{1,\dots, n_0-1\} \big)$ in the block form
{\fontsize{9}{10}
\begin{equation}
\label{eq:lDecompoJ0}
\begin{pmatrix}
S_0 & \bar{X}_0\\
X_0 & Q_0
\end{pmatrix}=
\renewcommand*{\arraystretch}{1.3}
 \left(
 \begin{array}{cc|cccc}
b(0)  & 0   & a(1) &  0 & \dots & 0\\
0 & b(0)  & 0 & 0 &  \dots & a(1)\\
 \hline
a(n_0) & 0 &  b(1) & a(2) &   \dots & 0  \\
0 & 0 &  a(n_0-1) & b(2) &    \ddots &  \vdots  \\
\vdots & \vdots  &  \vdots & \ddots & \ddots &   a(n_0-1)  \\
0 & a(n_0)  & 0 & 0 & a(2) & b(1) 
 \end{array}
 \right). 
\end{equation}}%
\\
We observe that $S_0$ is multiple of the identity matrix $b(0)I$ and that $Q_0=\jacobi^D_{cs}$. Let $z \in \rho(\jacobi^D_{cs})$, an element of the resolvent set of $\jacobi^D_{cs}$, then the Schur complement of $\jacobi_{cs}$ with respect to the decomposition (\ref{eq:lDecompoJ0}) is given by 
\begin{equation}\label{schurCompForJ0}
    Schur_{\ell(\partial G_{cs})}(\jacobi_{cs}):=zI - S_0 -  \bar{X}_0(zI-Q_0)^{-1} X_0.
\end{equation}
Our first result states that the Schur complement of $\jacobi_{cs}$ can be expressed using the polynomials  $P_{n_0}(z)$ and $P^D_{n_0-1}(z)$.
\begin{theorem}
\label{thm:SchurToPoly}
For each  $z \in \rho(\jacobi^D_{cs})$ we have  
%an element of the resolvent set of $\jacobi^D_{cs}$, then the Schur complement of $\jacobi_{cs}$ given in (\ref{schurCompForJ0}) satisfies the identity,
\begin{equation}\label{schurCompIdentity}
    Schur_{\ell(\partial G_{cs})}(\jacobi_{cs})=\phi_{\jacobi_{cs}}(z)\Big(  R_{\jacobi_{cs}}(z)I -\Delta_{0} \Big),
\end{equation}
where the functions $R_{\jacobi_{cs}}(z)$ and $\phi_{\jacobi_{cs}}(z)$ are  given by
\begin{align}
\label{eq:SpecDeciFunc}
R_{\jacobi_{cs}}(z):=1- \frac{P_{n_0}(z) }{\prod^{n_0}_{i=1}a(i)}, \quad \quad \phi_{\jacobi_{cs}}(z):= - \frac{\prod^{n_0}_{i=1}a(i)}{P^D_{n_0-1}(z)}.
\end{align}
\end{theorem}
\begin{remark}
Note that  $P^D_{n_0-1}(z) = \textbf{det} (zI-\jacobi^D_{cs})$. In particular, the domain of $\phi_{\jacobi_{cs}}(z)$ is given by $\rho(\jacobi^D_{cs})$.
\end{remark}
Before proving this result, we first establish an auxiliary Lemma.
\begin{lemma}
\label{lem: computingSchurForJ0}
Let $\bar{X}_0$, $Q_0$, and $X_0$ be the block matrices defined in (\ref{eq:lDecompoJ0}) and $z \in \rho(\jacobi^D_{cs})$. The following identity holds,
\begin{align}
\bar{X}_0 \big( zI-Q_0 \big)^{-1} X_0  =- \frac{\prod^{n_0}_{i=1}a(i)}{P^D_{n_0-1}(z)}
\begin{pmatrix}
-\frac{P^D_{n_0-2}(z)}{\prod^{n_0-1}_{i=2}a(i)}  & -1 \    \\
 \ -1   & -\frac{P^D_{n_0-2}(z)}{\prod^{n_0-1}_{i=2}a(i)} 
\end{pmatrix}.
\end{align}
In particular, this implies $Schur_{\ell(\partial G_{cs})}(\jacobi_{cs}) \in \mathcal{CSJ}$ for all $z \in \rho(\jacobi^D_{cs})$.
\end{lemma}
\begin{proof}
Recall that $Q_0$ is the Jacobi matrix $\jacobi^D_{cs}$. We invert the matrix $(zI-\jacobi^D_{cs})$ using Cramer's rule and restrict the computations to the relevant entries in the adjugate matrix, i.e. $\big( zI-Q_0 \big)^{-1}$ is equal

{\fontsize{9}{10} 
\begin{align}
\frac{1}{P^D_{n_0-1}(z)}
\begin{pmatrix}
P^D_{n_0-2}(z) & \ast & \dots &   \ast & (-1)^{n_0}(-1)^{n_0-2}\prod^{n_0-1}_{i=2}a(i)   \\
\ast & \ast & \dots &   \ast & \ast  \\
\vdots & \vdots & \vdots &    \vdots &  \vdots  \\
\ast & \ast & \dots & \ast &   \ast  \\
(-1)^{n_0}(-1)^{n_0-2}\prod^{n_0-1}_{i=2}a(i)  & \ast & \dots & \ast & P^D_{n_0-2}(z)
\end{pmatrix}.
\end{align}}%
For $X_0$ and $\bar{X}_0$ we have
\begin{equation}
X_0 =
\renewcommand*{\arraystretch}{1.3}
 \left(
 \begin{array}{cc}
a(n_0) & 0   \\
0 & 0   \\
\vdots & \vdots    \\
0 & a(n_0)   
 \end{array}
 \right), \quad \quad
  \bar{X}_0=
\renewcommand*{\arraystretch}{1.3}
 \left(
 \begin{array}{cccc}
 a(1) &  0 & \dots & 0\\
 0 & 0 &  \dots & a(1)
 \end{array}
 \right).
\end{equation}
A direct computation gives
\begin{align}
\bar{X}_0 \big( zI-Q_0 \big)^{-1} X_0  = \frac{a(1)a(n_0)}{P^D_{n_0-1}(z)}
\begin{pmatrix}
P^D_{n_0-2}(z)  & \prod^{n_0-1}_{i=2}a(i) \    \\
 \ \prod^{n_0-1}_{i=2}a(i)   & P^D_{n_0-2}(z) 
\end{pmatrix}.
\end{align}
\end{proof}

We can now we prove Theorem \ref{thm:SchurToPoly}.
\begin{proof}[Proof of Theorem \ref{thm:SchurToPoly}]
By Lemma \ref{lem: computingSchurForJ0} and (\ref{eq:SpecDeciFunc}), we have
\begin{align*}
\bar{X}_0 \big( zI-Q_0 \big)^{-1} X_0  &=\phi_{\jacobi_{cs}}(z)
\begin{pmatrix}
-\frac{P^D_{n_0-2}(z)}{\prod^{n_0-1}_{i=2}a(i)}  & -1 \    \\
 \ -1   & -\frac{P^D_{n_0-2}(z)}{\prod^{n_0-1}_{i=2}a(i)} 
\end{pmatrix}
\\
&  =\phi_{\jacobi_{cs}}(z)
\begin{pmatrix}
1  & -1 \    \\
 \ -1   &1
\end{pmatrix}
-\phi_{\jacobi_{cs}}(z)\big(1 +\frac{P^D_{n_0-2}(z)}{\prod^{n_0-1}_{i=2}a(i)}\big)
\begin{pmatrix}
1 & 0  \\
0 & 1
\end{pmatrix}.
\end{align*}
Recall that $S_0$ is multiple of the identity matrix $b(0)I$, hence

\begin{align*}
 Schur_{\ell(\partial G_{cs})}(\jacobi_{cs})  &=
 zI - S_0 -  \bar{X}_0(zI-Q_0)^{-1} X_0
\\
&  =
\phi_{\jacobi_{cs}}(z)\Big(\frac{z -b(0)}{\phi_{\jacobi_{cs}}(z)}  +1 +\frac{P^D_{n_0-2}(z)}{\prod^{n_0-1}_{i=2}a(i)}\Big)
I
-\phi_{\jacobi_{cs}}(z)
\Delta_{0}.
\end{align*}
Using (\ref{eq:SpecDeciFunc}) and the identity (\ref{eq:identityComputDet1}), we compute

\begin{align*}
\frac{z -b(0)}{\phi_{\jacobi_{cs}}(z)}  +1 +\frac{P^D_{n_0-2}(z)}{\prod^{n_0-1}_{i=2}a(i)}& =  -\frac{\big(z -b(0)\big)P^D_{n_0-1}(z)}{\prod^{n_0}_{i=1}a(i)} +1+ \frac{a(1)a(n_0)P^D_{n_0-2}(z)}{\prod^{n_0}_{i=1}a(i)} \\
&=1 -  \frac{ P_{n_0}(z)}{\prod^{n_0}_{i=1}a(i)}  = R_{\jacobi_{cs}}(z)
\end{align*}
\end{proof}

\subsection{Main technical tools}\label{sec:SpectDeciFunctsection}
This section is devoted to analyzing the polynomial $R_{\jacobi_{cs}}(z)$. In Section \ref{sec:DeciSection} we will see that $R_{\jacobi_{cs}}(z)$ is a \textit{spectral decimation function}, which represents a significant tool in the spectral analysis on fractals and self-similar graphs. The following theorem is the main result of this section and summarizes the relevant features of $R_{\jacobi_{cs}}(z)$ when investigating spectral properties of Jacobi matrices on graphs.
\begin{theorem}
\label{thm:propertiesOfSpecDeci}
The following statements hold:
\begin{enumerate}
\item[(1)] For all $z \in \complex$, we have 
\begin{align*}
 \text{det} \big(zI-\jacobi_{cs} )  \text{det}(zI-\jacobi^D_{cs}) = ( \prod^{n_0}_{i=1}a(i)  )^2R_{\jacobi_{cs}}(z)  \big(R_{\jacobi_{cs}}(z)- 2  \big) .
 \end{align*}
  In particular, $z \in \sigma(\jacobi_{cs}) \cup \sigma(\jacobi^D_{cs})$ if and only if $R_{\jacobi_{cs}}(z) \in \{0,2\}$.
\item[(2)] If $z \in \sigma(\jacobi_{cs}) \cap \sigma(\jacobi^D_{cs})$, then $z$ is a critical point of $R_{\jacobi_{cs}}$, i.e. $R'_{\jacobi_{cs}}(z)=0$.
\item[(3)] Let $z$ be a critical point of $R_{\jacobi_{cs}}(z)$, then $R_{\jacobi_{cs}}(z) \notin (0,2)$. In particular, there exist $n_0$ branches of the inverse $R^{-1}_{\jacobi_{cs}}$, that are defined and continuous in the domain $[0,2]$.
\end{enumerate}
\end{theorem}

To prove this result, we make the following elementary observation.

\begin{lemma}
\label{lem:criticalValuesOfPn0}
The polynomial $P_{n_0}(z)$ has $n_0$ simple real roots. The critical points of $P_{n_0}(z)$ are also real and interlaced between its roots. Moreover, if $z_c$ is a critical point of $P_{n_0}(z)$, then we have
\begin{align}
P_{n_0}(z_c) P^{''}_{n_0}(z_c) < 0.
\end{align}
\end{lemma}

We now have all the tools needed to prove the main result of this section.
\begin{proof}[Proof of Theorem \ref{thm:propertiesOfSpecDeci}] (1) Use Lemma \ref{Lem:PhiNenner} and the following computation {\fontsize{9}{10} 
\begin{align*}
\frac{1}{\big( \prod^{n_0}_{i=1}a(i)  \big)^2} \big[\big(P_{n_0}(z)\big)^2-\big( \prod^{n_0}_{i=1}a(i)  \big)^2 \big]&=\Big(\frac{P_{n_0}(z)}{ \prod^{n_0}_{i=1}a(i)} - 1  \Big) \Big(\frac{P_{n_0}(z)}{ \prod^{n_0}_{i=1}a(i)} + 1  \Big) \\
&= R_{\jacobi_{cs}}(z)  \big(R_{\jacobi_{cs}}(z)- 2  \big),
\end{align*}}%
where in the last equality, we used the definition of the $R_{\jacobi_{cs}}$ given in (\ref{eq:SpecDeciFunc}). 
\begin{enumerate}
\item[(2)] The second statement follows from differentiating the identity in part (1):
%, i.e. $2 R_{\jacobi_{cs}}(z) R'_{\jacobi_{cs}}(z) - 2 R'_{\jacobi_c}(z)$ is equal to 
{\fontsize{9}{10}
\begin{align}
%\label{eq:differentOfdeterminantIden}
2 R_{\jacobi_{cs}}(z) R'_{\jacobi_{cs}}(z) - 2 R'_{\jacobi_c}(z)= \notag \\
  \frac{ [\text{det} \big(zI-\jacobi_{cs} )]^{'}  \text{det}(zI-\jacobi^D_{cs}) +  \text{det} \big(zI-\jacobi_{cs} )  [\text{det}(zI-\jacobi^D_{cs})]^{'} }{\big( \prod^{n_0}_{i=1}a(i)  \big)^2}
\label{eq:differentOfdeterminantIden}
\end{align} }%
\item[(3)] If the critical point $z_{c}$ is in $ \sigma(\jacobi_{cs}) \cup \sigma(\jacobi^D_{cs})$, then the first part of the statement holds by part (1). We assume $z_{c} \notin \sigma(\jacobi_{cs}) \cup \sigma(\jacobi^D_{cs})$ and define
\begin{align}
D(z) := \frac{\text{det} \big(zI-\jacobi_{cs} )  \text{det}(zI-\jacobi^D_{cs})}{\big( \prod^{n_0}_{i=1}a(i)  \big)^2}
\end{align}
Note that the polynomial $D(z)$ is of degree $2n_0$ and the roots of $D(z)$  are given by the eigenvalues  $ \sigma(\jacobi_{cs}) \cup \sigma(\jacobi^D_{cs})$. By (\ref{eq:differentOfdeterminantIden}), we see that $z_{c}$ is also a critical point of $D(z)$. Moreover, $z_{c}$ lies between two roots of $D(z)$ due to the assumption $z_{c} \notin \sigma(\jacobi_{cs}) \cup \sigma(\jacobi^D_{cs})$. Using Lemma \ref{Lem:PhiNenner}, we compute
\begin{align}
D^{''}(z) = \frac{2 \big(P^{'}_{n_0}(z)\big)^2+2 P_{n_0}(z) P^{''}_{n_0}(z)}{\big( \prod^{n_0}_{i=1}a(i)  \big)^2}
\end{align}
By (\ref{eq:SpecDeciFunc}), we have $P^{'}_{n_0}(z_{c})=0$ and Lemma \ref{lem:criticalValuesOfPn0} implies that $z_{c}$ is a local maximum of $D(z)$, i.e.
\begin{align}
D^{''}(z_{c}) = \frac{2 P_{n_0}(z_{c}) P^{''}_{n_0}(z_{c})}{\big( \prod^{n_0}_{i=1}a(i)  \big)^2}<0
\end{align}
In particular, we have $D(z_{c})>0$. We rewrite the identity in part (1),
\begin{align}
(R_{\jacobi_{cs}}(z_{c}) -1)^2 = R_{\jacobi_{cs}}(z_{c})^2 - 2 R_{\jacobi_{cs}}(z_{c}) + 1 = D(z_{c}) + 1,
\end{align}
which gives $R_{\jacobi_{cs}}(z_{c})   = 1 \pm \sqrt{D(z_{c}) + 1} \  \notin (0,2)$.
\end{enumerate}
\end{proof}

\section{Spectral decimation of piecewise centrosymmetric Jacobi matrices on graphs}
\label{sec:DeciSection}

We start by recalling some facts about of the spectral decimation method following the framework set forth in ~\cite{MalozemovTeplyaev2003}, as it is the best suited for our results.  Let $\mathcal{H}$ and $\mathcal{H}_0$ be Hilbert spaces, and $U:\mathcal{H}_0 \to \mathcal{H}$ be an isometry. Suppose  $H$ and $H_0$ are bounded linear operators on  $\mathcal{H}$ and $\mathcal{H}_0$, respectively, and that $\phi,\psi$ are complex-valued functions. We call the operator $H$ \textit{spectrally similar} to the operator $H_0$ with functions $\phi$ and $\psi$ if \cite[Definition 2.1]{MalozemovTeplyaev2003}
\begin{equation}
\label{eq:OriginalSpectSimi}
U^{\ast}(H-zI)^{-1}U=(\phi(z)H_0 - \psi(z)I)^{-1},
\end{equation}
for all $z \in \complex$ such that the two sides of~\eqref{eq:OriginalSpectSimi} are well defined. Note, in particular, that for $z$ in the domain of both $\phi$ and $\psi$ and satisfying $\phi(z)\neq0$ we have $z\in\rho(H)$ (the resolvent set of $H$) if and only if $R(z)=\frac{\psi(z)}{\phi(z)}\in\rho(H_0)$.  We call $R(z)$ the {\em spectral decimation function}. The functions $\phi(z)$ and $\psi(z)$ are difficult to determine directly from the structure of the considered fractal or graph, but they can be computed effectively using a Schur complement (several examples may be found in~\cite{BajorinVibrationSpectra2008, BajorinVibration3Ngasket2008}). Identifying $\mathcal{H}_0$  with a closed subspace of $\mathcal{H}$ via $U$, let $\mathcal{H}_1$ be the orthogonal complement and decompose $H$ on   $\mathcal{H}=\mathcal{H}_0 \oplus \mathcal{H}_1$ in the block form
\begin{equation}
\label{eq:BlockDecompo}
H=\begin{pmatrix}
S & \bar{X}\\
X & Q
\end{pmatrix}.
\end{equation}

\begin{lemma}(\cite[Lemma 3.3]{MalozemovTeplyaev2003})
\label{lem:Lemma33}
For $z \in \rho(H)\cap \rho(Q)$ the operators $H$ and $H_0$ are spectrally similar if and only if
the Schur complement of $H$ (with respect to the decomposition (\ref{eq:BlockDecompo})), given by $Schur_{\mathcal{H}_0}(H):=zI-S -   \bar{X} (zI-Q)^{-1}X$, satisfies
\begin{equation}
\label{eq:schurComp1}
    Schur_{\mathcal{H}_0}(H)=   \psi(z)I-\phi(z) H_0.
\end{equation}
\end{lemma}
The following set plays a crucial role in the spectral decimation method. 
\begin{definition}
\label{exceptionalSet}
The \textit{exceptional set} of $H$ is given by $\mathscr{E}_{H}:=\{z \in \complex \ \rvert \ z \in \sigma(Q) \text{ or } \phi(z)=0 \}$
\end{definition}
The importance of the concept of \textit{spectral similarity} stems from the following result.
\begin{proposition}(\cite[Theorem 3.6]{MalozemovTeplyaev2003})
\label{thm:mainMalozemovTepl}
Let $H$ be spectrally similar to $H_0$ with functions $\phi$ and $\psi$ and $z \notin \mathscr{E}_{H}$. Then
\begin{enumerate}
\item $R(z) \in \rho(H_0)$ if and only if $z \in \rho(H)$.
\item $R(z)$ is an eigenvalue of $H_0$ if and only if $z$ is an eigenvalue of $H$. Moreover, there is a one-to-one map
\begin{align}
f_0 \mapsto f := f_0 - (zI-Q)^{-1}X f_0
\end{align}
from the eigenspace of $H_0$ corresponding to $R(z)$ onto the eigenspace of $H$ corresponding to $z$.
\end{enumerate}
\end{proposition}

We now apply this framework of the spectral decimation method to the substitution graph $G$ and operator $\jacobi = F_{G_{p}}(\Delta_{\mathbf{p}}, \jacobi_{cs})$ obtained from a graph $G_p$ equipped with a probabilistic Laplacian $\Delta_p$, and a graph  $G_{cs}$ be a finite path graph associated to a centrosymmetric weight matrix $\jacobi_{cs}$. 

%. To this end, we suppose $G_{p}$ is defined as in section \ref{sec:ProbabilisticgraphLaplacian} and fixed. Let $G_{cs}$ be a finite path graph associated to a centrosymmetric weight matrix $\jacobi_{cs}$ and $G=(V(G),E(G))$ is the resulting graph, when substituting $G_{cs}$ in $G_{p}$ as described in section \ref{sec:WeightedSubstitutionOperator}. Let $\Delta_{\mathbf{p}} \in \mathcal{L}_{G_{p}}$ and the weight matrix associated with the graph $G$ is given by $\jacobi = F_{G_{p}}(\Delta_{\mathbf{p}}, \jacobi_{cs})$.
\begin{assumption}
\label{assump:HilbertspacesAndIsometry}
We assume:
\begin{enumerate}
\item There exists a Hilbert space of  $\complex$-valued functions on $V(G_{p})$, which we denote by  $\ell^2( G_{p}, d \pi_{p})$, such that the probabilistic graph Laplacian $\Delta_{\mathbf{p}}$ is bounded and self-adjoint.
\item There exists a Hilbert space of  $\complex$-valued functions on $V(G)$, which we denote by $\ell^2( G, d \pi_{J})$, such that $\jacobi $ is bounded and self-adjoint.
\item There exists an isometry $U: \ell^2( G_{p}, d \pi_{p}) \to \ell^2( G, d \pi_{J})$.
\end{enumerate}
\end{assumption}
\begin{remark}

\begin{enumerate}
\item
The self-adjointness assumption is by no means necessary. In fact, the results in \cite{MalozemovTeplyaev2003} 
%we are referring to 
hold for more general operators. Section \ref{sec:OnedimPathGraph} shows how to construct such Hilbert spaces in the one-dimensional path graphs case, for which Assumption \ref{assump:HilbertspacesAndIsometry} holds. The key idea is to equip the set of vertices $V(G_{p})$ and $V(G)$ with a measure satisfying a Kolmogorov’s cycle type condition, see~\eqref{eq:measureDef} for more details.
\item
\label{Rem:specSubsetOfInterval}
We note that Assumption \ref{assump:HilbertspacesAndIsometry} combined with Perron-Frobenius theorem implies that $\sigma(\Delta_{\mathbf{p}}) \subset [0,2]$, as $\Delta_{\mathbf{p}}$ is a stochastic matrix (with spectrum in $[-1,1]$) shifted by the identity, see also \cite[Remark 5.9]{MalozemovTeplyaev2003}.
\end{enumerate}
\end{remark}

\begin{proposition}
\label{prop:SpectralSimilarity}
The following statements hold:
\begin{enumerate}
\item[(1)] For the exceptional set, we have $\mathscr{E}_{\jacobi_{cs}}= \sigma \big(\jacobi^D_{cs} \big)$.
\item[(2)]  $\jacobi_{cs}$ is spectrally similar to $ \Delta_{0}$ (given in (\ref{eq:TrivialProbLap})) with the functions $\phi_{\jacobi_{cs}}(z)$ defined in (\ref{eq:SpecDeciFunc}), and $\psi_{\jacobi_{cs}}(z):=\phi_{\jacobi_{cs}}(z)R_{\jacobi_{cs}}(z)$, where $R_{\jacobi_{cs}}(z)$ is given by (\ref{eq:SpecDeciFunc}). In particular,  $R_{\jacobi_{cs}}(z)$ is the corresponding spectral decimation function.
\item[(3)] $\jacobi = F_{G_{p}}(\Delta_{\mathbf{p}}, \jacobi_{cs})$ is spectrally similar to $\Delta_{\mathbf{p}}$ with the same functions $\phi_{\jacobi_{cs}}(z)$ and  $\psi_{\jacobi_{cs}}(z)$ as in part (2). The associated spectral decimation function and exceptional set coincide with $R_{\jacobi_{cs}}(z)$ and $\mathscr{E}_{\jacobi_{cs}}$, respectively.
\end{enumerate}
\end{proposition}

\begin{proof} (1) This statement follows by Definition \ref{exceptionalSet} and $\phi_{\jacobi_{cs}}(z) \neq 0$ (see (\ref{eq:SpecDeciFunc}) and recall the assumption $a(i)\neq 0$). 
\begin{enumerate}
\item[(2)] This assertion follows by Lemma \ref{lem:Lemma33} and Theorem \ref{thm:SchurToPoly}.
\item[(3)]  This final statement   is an immediate consequence of Lemma \cite[Lemma 3.10]{MalozemovTeplyaev2003} and its main idea can be sketched as follows. 
Using the isometry $U$, we decompose $\jacobi$ on   $ \ell^2( G, d \pi_{J})= U\big( \ell^2( G_{p}, d \pi_{p}) \big) \oplus U\big( \ell^2( G_{p}, d \pi_{p}) \big)^{\bot}$ in the block form (\ref{eq:BlockDecompo}) and denote the Schur complement of $\jacobi$ with respect to the decomposition (\ref{eq:BlockDecompo}) by  $Schur_{ \ell ( V(G_{p}) )}(\jacobi):=zI - S -  \bar{X}(zI-Q)^{-1} X$. We prove that the Schur complement is preserved under the substitution operator $F_{G_{p}}$, in the sense
\begin{align}
F_{G_{p}}\big(\Delta_{\mathbf{p}},  Schur_{\ell(\partial G_{cs})}(\jacobi_{cs}) \big) = Schur_{ \ell ( V(G_{p}) )}(\jacobi).
\end{align}
Proposition \ref{prop:PropertiesOfFstar} and Theorem \ref{thm:SchurToPoly} imply
\begin{align*}
Schur_{ \ell ( V(G_{p}) )}(\jacobi) &= \phi_{\jacobi_{cs}}(z) \ F_{G_{p}}\Big(\Delta_{\mathbf{p}},   R_{\jacobi_{cs}}(z)I -\Delta_{0} \Big) \\
&= \phi_{\jacobi_{cs}}(z)\Big(     R_{\jacobi_{cs}}(z) \     F_{G_{p}}\big(\Delta_{\mathbf{p}},   I  \big)  - F_{G_{p}}\big(\Delta_{\mathbf{p}},  \Delta_{0} \big)  \Big) \\
&= \phi_{\jacobi_{cs}}(z)\Big(     R_{\jacobi_{cs}}(z) I  - \Delta_{\mathbf{p}}   \Big).
\end{align*}
The second statement follows then by Lemma \ref{lem:Lemma33}.
\end{enumerate}
\end{proof}
The main result of this section is the following.
\begin{theorem}
\label{coro:RenoEquAndSpectra}
The following statements hold:
\begin{enumerate}
\item[(1)] 
\label{key1}
$R^{-1}_{\jacobi_{cs}} \big(  \sigma(\Delta_{\mathbf{p}}) \setminus \{0,2\} \big) 
\subset
\sigma \big( \jacobi  \big)
\subset
R^{-1}_{\jacobi_{cs}} \big(  \sigma(\Delta_{\mathbf{p}}) \cup \{0,2\} \big) $. 
\item[(2)] 
\label{key2}
The resolvent operators satisfy the following renormalization identity
\begin{align}
\label{eq:RenoEquation}
U^{\ast}\big(zI - \jacobi\big)^{-1}U=-\frac{P^D_{n_0-1}(z)}{\prod^{n_0}_{i=1}a(i)}\big(R_{\jacobi_{cs}}(z) I- \Delta_{\mathbf{p}}\big)^{-1}
\end{align}
for $z \in \rho(\jacobi)$ and $z \notin \mathscr{E}_{\jacobi_{cs}}= \sigma \big(\jacobi^D_{cs} \big)$.
\item[(3)] 
\label{key3}
%{\color{blue} 
Polynomial $R_{\jacobi_{cs}}$ gives an isomorphism of the types of the spectral measures of the operators $\jacobi$ and 
$\Delta_{\mathbf{p}}$	outside of the finite exceptional set $R^{-1}_{\jacobi_{cs}} \big(    \{0,2\} \big)$.
% }
\end{enumerate}
\end{theorem}

\begin{proof}
Theorem \ref{thm:propertiesOfSpecDeci} asserts that is $ \mathscr{E}_{\jacobi_{cs}}  \subset R^{-1}_{\jacobi_{cs}}\{0,2\}$. 
Statements (1) and (2) follow by Proposition \ref{thm:mainMalozemovTepl} and Proposition \ref{prop:SpectralSimilarity}. 

%{\color{blue} 
Statement (3) is proved following the same lines as \cite[Theorem 2.3]{StrichartzTeplyaev2012}.  We use \eqref{eq:RenoEquation}, \eqref{eq:OriginalSpectSimi} and 
	the Schur complement in Lemma \ref{lem:Lemma33} with the block decomposition 
	\eqref{eq:BlockDecompo}.  Using the standard general theory \cite[Section VIII.7]{ReedSimonIFunctionalAnalysis} or \cite[Chapter 3]{TeschlQM2014}, 	
	for a self-adjoint operator $H$ on a Hilbert space $\mathcal{H}$ the spectral measure $\mu_{\psi}$ of $\psi \in \mathcal{H}$ is uniquely defined by the Herglotz function $F_{\psi}(z)$, a Borel transform of a finite Borel measure, 
$
			F_{\psi}(z) = \bra{\psi}\ket{(H-z)^{-1}  \psi}_{\mathcal{H}} = \int_{\rr} \frac{1}{\lambda - z} d \mu_{\psi}(\lambda)
$, $z\in\cc\setminus\rr$. The measure $\mu_{\psi}$ is unique  by the Stieltjes inversion formula, which implies the result. 
	In particular, if 	$\ast \ast \in \{pp,sc,ac\}$ stands for   \emph{pure point, singular continuous, absolutely continuous} spectrum, then 
for  $z \notin \mathscr{E}_{\jacobi_{cs}}$   we have $R(z) \in \sigma_{\ast \ast}(\Delta_{\mathbf{p}}) \text{ if and only if } z \in \sigma_{\ast \ast}(\jacobi)$.
%}
\end{proof}

\begin{remark}
Theorem 
\ref{coro:RenoEquAndSpectra} can be extended to a formula for the spectral projections and the spectral resolution of the identity, similarly to \cite[Theorem 2.3]{StrichartzTeplyaev2012}, using the general theory of Schur complement for self-adjoint operators, see \cite{ReedSimonIFunctionalAnalysis,BajorinVibration3Ngasket2008,Teplyaev1998,KZ19,LenzT,Str89,Str92,revisited}. Useful resolvent estimates are given in \cite{res1,res2}.
\end{remark}

\begin{corollary}
\label{IninitelySubstitutionsTHM}
If the points $0$ and $2$ are limit points of $\sigma(\Delta_{\mathbf{p}})$, then $\sigma \big( \jacobi  \big) = R^{-1}_{\jacobi_{cs}} \big(  \sigma(\Delta_{\mathbf{p}}) \big)$.
\end{corollary}
\begin{proof}
Theorem  \ref{thm:propertiesOfSpecDeci} asserts the existence and continuity of the branches of the inverse $R^{-1}_{\jacobi_{cs}}$ on the interval $[0,2]$. We conclude {\fontsize{8.5}{10}
\begin{align*}
R^{-1}_{\jacobi_{cs}} \big(  \sigma(\Delta_{\mathbf{p}}) \big)  =
\overline{ R^{-1}_{\jacobi_{cs}} \big(  \sigma(\Delta_{\mathbf{p}}) \setminus \{0,2\} \big) }
\subset \sigma \big( \jacobi  \big) 
\subset \overline{ R^{-1}_{\jacobi_{cs}} \big(  \sigma(\Delta_{\mathbf{p}})  \cup \{0,2\} \big) }
= R^{-1}_{\jacobi_{cs}} \big(  \sigma(\Delta_{\mathbf{p}}) \big).
\end{align*}}%
\end{proof}

%{\fontsize{9}{10} }%

 \section{The one-dimensional path graphs case}\label{sec:OnedimPathGraph}
 
 %{\color{blue}
 
\begin{figure}[htp]
\centering
%\hspace*{-4cm} 
 \includegraphics[width=1.\textwidth]{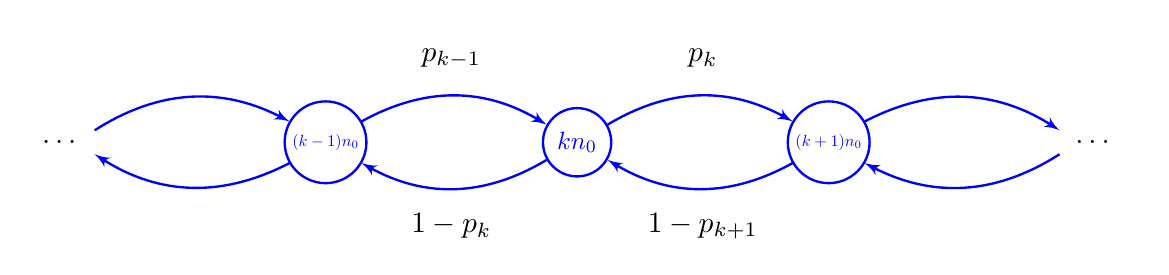}
\caption{ The directed path graph $G_{p} = (V(G_{p}),E(G_{p}))$.}
\label{fig:PartitionOfGZero}
\end{figure}

 In this section, we illustrate the results of the previous sections by focusing on  one-dimensional path graphs. We apply the framework of Section \ref{sec:DeciSection} and start by addressing the question concerning the existence of the Hilbert spaces and the isometry mapping in Assumption \ref{assump:HilbertspacesAndIsometry}. Let $G_{cs}=(V(G_{cs}),E(G_{cs}))$ be defined as in Section  \ref{sec:centrosymmetryJacobis}. Suppose that $G_{cs}$ consists of $n_0+1$ vertices for fixed $n_0 \geq 1$ and the associated weight matrix $\jacobi_{cs} \in \mathcal{CSJ}$ is of the form (\ref{eq:BasicJacobi}). Next, we restrict the definition of $G_{p} = (V(G_{p}),E(G_{p}))$ as given in Section \ref{sec:ProbabilisticgraphLaplacian} to path graphs. We distinguish three cases, when $G_{p} $ is a finite, semi-infinite, or infinite path graph. Accordingly, the set of vertices is given by % to which case is considered, we set for the set of vertices
\begin{equation}
%	\label{eq:setOfVertices}
	V(G_{p} )=\{0,n_0, \dots, k_0 n_0\},  \quad V(G_{p} )= n_0\integers_{+},\text{ or }   \  \   V(G_{p} )= n_0\integers.
\end{equation}
The set of edges $E(G_{p} )$ is then given by $\big\{(kn_0,kn_0+n_0) \ \rvert \ kn_0,(k+1)n_0 \in  V(G_{p} ) \big\} \cup \big\{(kn_0+n_0,kn_0) \  \rvert \ kn_0,(k+1)n_0 \in  V(G_{p} )  \big\}$.
We simplify the notation and set for the transition probabilities $p_k : = p(kn_0,kn_0+n_0)$. In particular, the third equation in conditions (\ref{eq:CodForTransProb}) implies that $p(kn_0,kn_0-n_0)=1-p_k$, see Figure \ref{fig:PartitionOfGZero}. 
\begin{figure}[htp]
\centering
%\hspace*{-4cm} 
 \includegraphics[width=1.\textwidth]{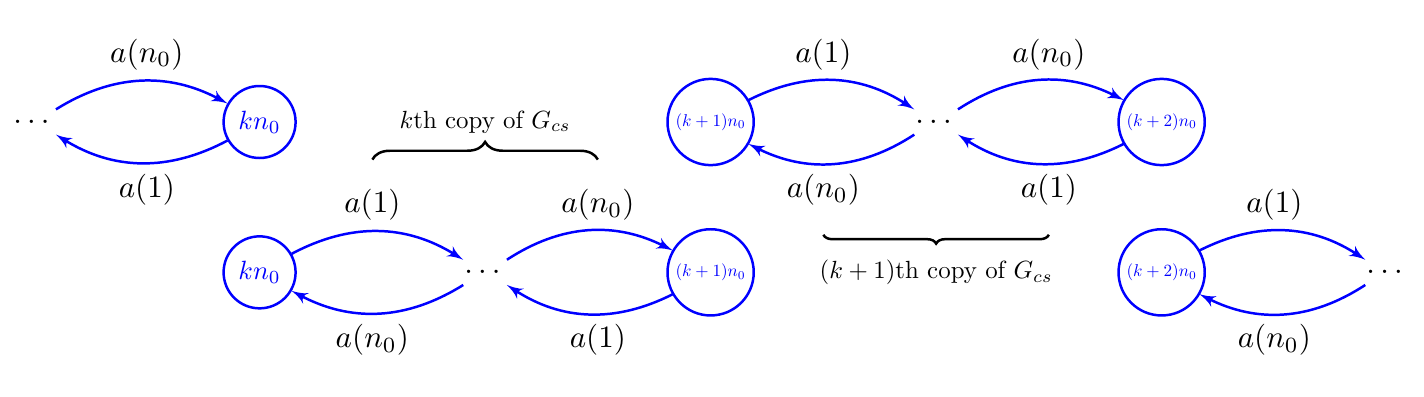}
\caption{Constructing the graph $G$ by substituting copies of $G_{cs}$ between adjacent vertices in $G_{p}$.}
\label{fig:PartitionOfG}
\end{figure}
Let $G=(V(G),E(G))$ be the constructed graph when substituting copies of $G_{cs}$ between adjacent vertices in $G_{p}$ as explained in Section~\ref{sec:WeightedSubstitutionOperator}.
Depending on  $G_{p}$, we distinguish three cases, where $G$ is a finite, semi-infinite or infinite path graph. Accordingly, the set of vertices is given by %we set for the set of vertices
	\begin{equation}
%	\label{eq:setOfVertices}
	V(G)=[0,n] \cap  \integers_{+},  \quad V(G)= \integers_{+},\text{ or }   \  \   V(G)= \integers,
	\end{equation}
where $n \in \nn$, $n \geq 1$. In the finite graph case we assume that $n_0$ is a divisor of $n$, i.e. $n=k_0 n_0$, and $k_0$ is the number of the copies of $G_{cs}$ in $G$.  We regard the vertices $\{kn_0, \dots , (k+1)n_0\}$ as the vertices of the $k$th copy of  $G_{cs}$, see Figure \ref{fig:PartitionOfG}.  The weight matrix associated with the graph $G$ is a tridiagonal matrix, that is created by the substitution operator $\jacobi =F_{G_{p}}(\Delta_{\mathbf{p}}, \jacobi_{cs}) $. The operators $\Delta_{\mathbf{p}}$ and $\jacobi$ are assumed to act on the following Hilbert spaces,{\fontsize{8}{10}
\begin{align*}
\begin{cases}
 \   \ell^2( G_{p}, d \pi_{p}) := \Big\{ f:  V(G_{p}) \to \complex \ \big\rvert \ \sum_{x \in  V(G_{p})} \ \rvert f (x)\rvert^2 \pi_{p}(x) < \infty  \Big\}, \\ 
 \  \langle f, g \rangle_{p} = \sum_{x \in  V(G_{p})} \overline{f(x)} g(x)\pi_{p}(x), 
    \end{cases}
    \\
    \begin{cases}
 \   \ell^2( G, d \pi_{J}) := \Big\{ f: V(G) \to \complex \ \rvert \ \sum_{x \in  V(G)} \ |f (x)|^2  \pi_{J}(x) < \infty  \Big\}, \\ 
 \ \langle f, g \rangle_{J} = \sum_{x \in  V(G)} \overline{f(x)} g(x) \pi_{J}(x),
    \end{cases}
\end{align*}}%
respectively, where the set of vertices $V(G_{p})$ and $V(G)$ are equipped with the following measures satisfying the Kolmogorov’s cycle condition \cite{DurrettBookProbability2019, GrigoryanIntroGraphs2018,KellerLenzBook2021,KellyBookReversibility2011},
{\fontsize{9}{10}
\begin{align}
\label{eq:measureDef}
\begin{cases}
	\pi_{p}(0)=1, \quad \pi_{p}(kn_0) = \pi_{p}  \big( (k-1)n_0 \big) \frac{p_{k-1}}{p_k}, \  \  \  kn_0 \in V(G_{p})\backslash \{0\}, \\
           \pi_{J}(0)=1, \quad  \pi_{J}(x) = \pi_{J}(x-1) \frac{\jacobi(x-1,x)}{\jacobi(x,x-1)},    \ \  x \in V(G).
\end{cases}
\end{align} }%
\\
Note that if we set $\jacobi(x-1,x)=\jacobi(x,x-1)$ for all $x \in \integers$ in the infinite graph case, we obtain $\ell^2( G, d \pi_{J}) = \ell^2( \integers)$.
\begin{proposition}
\label{prop:SelfAdjointnessOfOps}
The operator $\jacobi$ (resp. $\Delta_{\mathbf{p}}$) is a bounded self-adjoint operator on $ \ell^2( G, d \pi_{J})$ (resp. $ \ell^2( G_{p}, d \pi_{p})$ ). Moreover, for $x \in V(G_{p}) \subset V(G)$, we have $\pi_{p}(x) = \pi_{J}(x)$.
\end{proposition}
\begin{proof}
Let  $W_n(f,g)  :=  \pi_{J}(n) \jacobi(n,n+1)    \big(    f(n)g(n+1)-f(n+1)g (n)  \big)$ be the $n$-th Wronskian of $f$ and $g$. Direct computation gives
\begin{equation*}
 \sum_{x=m}^n f(x) \jacobi g (x) \pi_{J}(x) -  \sum_{x=m}^n  \jacobi f(x)  g (x) \pi_{J}(x) = W_n(f,g)-W_{m-1}(f,g).
\end{equation*}
For $f,g \in \ell^2( G, d \pi_{J})$, we see that $ \langle f, \jacobi g \rangle_{J}-  \langle \jacobi f, g \rangle_{J} = \lim_{n \to \infty} W_n(f,g) - \lim_{m \to -\infty} W_{m-1}(f,g) =0$. Similar computations for the semi-infinite graph case and the proof is identical for $\Delta_{\mathbf{p}}$. The second statement follows by (\ref{eq:measureDef}) and the centrosymmetry  assumption on $\jacobi_{cs}$.
\end{proof}

Using Proposition \ref{prop:SelfAdjointnessOfOps}, we view  $\ell^2( G_{p}, d \pi_{p})$ as a subspace of $  \ell^2( G, d \pi_{J})$ via the identification 
%(justified by Proposition \ref{prop:SelfAdjointnessOfOps}, i.e $\pi_{p}(x) = \pi_{J}(x)$ for $x \in V(G_{p}) \subset V(G)$),
\begin{align*}
 \ell^2( G_{p}, d \pi_{p})   \    \cong \ \Big\{     \psi \in  \ell^2( G, d \pi_{J})  \ \rvert  \   \psi(x)=0 \ \forall x \in V(G)\backslash V(G_{p})       \Big\}.
\end{align*}
In this sense, $U:\ell^2( G_{p}, d \pi_{p}) \to \ell^2( G, d \pi_{J}) $ is the inclusion operator.
%\subsection{Finite path graphs}\label{sec:finiteGraphcase}
\begin{proposition}
\label{FinitelySubstitutionsTHM}
Assume $G_{p}$ is a finite path graph. Then $$\sigma \big( \jacobi  \big) = \big( R^{-1}_{\jacobi_{cs}} \big(  \sigma(\Delta_{\mathbf{p}})\big) \big\backslash \{0,2\} \big) \cup \sigma \big(\jacobi_{cs} \big)$$.
\end{proposition}
\begin{proof}
The constant vector $(1,1,1,\dots)^t$ is an eigenvector with the eigenvalue $0$ and the alternating vector $(1,-1, 1\dots)^T$ is an eigenvector with the eigenvalue $2$ for $\Delta_{\mathbf{p}}$. Hence, we have $\{0,2\} \subset \sigma(\Delta_{\mathbf{p}}) \subset [0,2]$. In particular, $\big \rvert \sigma(\Delta_{\mathbf{p}}) \backslash  \{0,2\} \big \rvert=k_0-1$. Note that $R_{\jacobi_{cs}}(z)$ is as $P_{n_0}(z)$ a polynomial of degree $n_0$. Theorem \ref{thm:propertiesOfSpecDeci} implies
\begin{align*}
	\Big \rvert R^{-1}_{\jacobi_{cs} } \Big( \sigma(\Delta_{\mathbf{p}})\backslash \{0,2\}  \Big) \Big \rvert =n_0(k_0-1) = n_0k_0 - n_0.
\end{align*}
Note that all the $n_0k_0 - n_0$ preimages are not in the exceptional set and therefore distinct eigenvalues of  $\jacobi$. Theorem \ref{thm:propertiesOfSpecDeci} asserts $ R^{-1}_{\jacobi_{cs} } \big(  \{0,2\}  \big)=\sigma(\jacobi_{cs}) \cup \sigma(\jacobi^D_{cs})$. By excluding the exceptional points   $\mathscr{E}_{\jacobi_{cs}}= \sigma \big(\jacobi^D_{cs} \big)$, we see that $R^{-1}_{\jacobi_0} \big(  \{0,2\} \big)$ generates $n_0+1$ distinct eigenvalues of  $\jacobi$, namely the eigenvalues in $\sigma \big(\jacobi_{cs}\big)$. We generated $(n_0k_0 - n_0 )+( n_0+1) = n_0k_0 +1$ distinct eigenvalues, which shows with a dimension argument that we completely determined the spectrum $\sigma\big( \jacobi  \big)$.
\end{proof}
%\begin{remark}
%comment and refer to the theorem in AMO papar
%\end{remark}
% \subsection{Constructing a sequence of self-similar probabilistic Laplacians}
%\label{sec:probabilisticLaplacians}
We now restrict our investigations to the case $G_{p}$ is a finite path graph and use the substitution operator $F_{G_{p}}$ to generate a multi-parameter family of self-similar probabilistic Laplacians. A key observation is that the centrosymmetric probabilistic Laplacians are invariant under the substitution operator in the sense that if $\Delta_{\mathbf{p}_1}, \Delta_{\mathbf{p}_2} \in  \mathcal{CSJ} $ with $\Delta_{\mathbf{p}_1} \in \mathcal{L}_{G_{p}}$  then $F_{G_{p}}(\Delta_{\mathbf{p}_1}, \Delta_{\mathbf{p}_2}) \in  \mathcal{CSJ} \cap \mathcal{L}_{G}$.
\begin{definition}
\label{def:selfSimiLapSeq}
Let $\Delta^{(1)}_{\mathbf{p}}  \in  \mathcal{CSJ} \cap \mathcal{L}_{G_{p}}$ be given. We define the sequence $\{\Delta^{(\ell)}_{\mathbf{p}}\}_{\ell \in \nn}  \subset  \mathcal{CSJ}$ inductively by
$\Delta^{(\ell+1)}_{\mathbf{p}}  := F_{G_{p}}(  \Delta^{(1)}_{\mathbf{p}}, \Delta^{(\ell)}_{\mathbf{p}} ) $ for $\ell \geq 1$.
We refer to $\{\Delta^{(\ell)}_{\mathbf{p}}\}_{\ell \in \nn}$ as a sequence of self-similar probabilistic Laplacians.
\end{definition}
If we initialize $\Delta^{(1)}_{\mathbf{p}}:=\Delta_{0}$, then $\Delta^{(\ell)}_{\mathbf{p}}$ is a constant sequence, i.e. $\Delta^{(\ell)}_{\mathbf{p}}=\Delta_{0}$ for all $\ell \geq 1$. For $k_0 \geq 2$ , $\ \mathbf{p}=\{p_0,p_1, \dots, p_{k_0}\}$ with $p_0=1$ and $p_{k_0}=0$, we initialize $\Delta^{(1)}_{\mathbf{p}}$  as the following $(k_0+1) \times (k_0+1)$-matrix{\fontsize{9}{10} 
\begin{align}
\label{eq:initialLapNandD}
\Delta^{(1)}_{\mathbf{p}} =
\begin{pmatrix}
1 & -1 & 0 &   \dots & 0   \\
p_1-1 & 1 & -p_1 &   \dots & 0  \\
0 & p_2-1 & 1 &    \ddots &  \vdots  \\
\vdots & \vdots & \ddots & \ddots &   p_1-1  \\
0 & 0 & 0 & -1 & 1 
\end{pmatrix},
\quad 
\Delta^{(1),D}_{\mathbf{p}} =
\begin{pmatrix}
 1 & -p_1 &   \dots & 0  \\
 p_2-1 & 1 &    \ddots &  \vdots  \\
 \vdots & \ddots & \ddots &   p_2-1  \\
 0 & 0 & -p_1 & 1  
\end{pmatrix}.
\end{align} }%
Note that the centrosymmetry assumption impose $ p_1 = 1-p_{k_0 -1}$, $ p_2 = 1-p_{k_0 -2}$, and so on.
Proposition \ref{prop:SpectralSimilarity} implies that $\Delta^{(1)}_{\mathbf{p}}$ is spectrally similar to $ \Delta_{0}$ with the functions $R_{\Delta^{(1)}_{\mathbf{p}}}$ and $\phi_{\Delta^{(1)}_{\mathbf{p}}}$ given by
\begin{align}
\label{eq:SpecDeciFuncForSeqOfLap}
R_{\Delta^{(1)}_{\mathbf{p}}}(z)=1- (-1)^{k_0}\frac{P_{k_0}(z) }{\prod^{k_0-1}_{i=0}p_i}, \quad \quad \phi_{\Delta^{(1)}_{\mathbf{p}}}(z)= (-1)^{k_0+1} \frac{\prod^{k_0-1}_{i=0}p_i}{P^D_{k_0-1}(z)},
\end{align}
where the polynomial $P_{k_0}$ (resp. $P^D_{k_0-1}$)
is defined by the recurrence relation (\ref{eq:recurrenceRelations}) corresponding to the matrix $\Delta^{(1)}_{\mathbf{p}} $ (resp. by the recurrence relation (\ref{eq:recurrenceRelationsDiri}) corresponding to the matrix $\Delta^{(1),D}_{\mathbf{p}} $). 

Note that in the special case of a symmetric random walk, i.e. $ p_k=\frac{1}{2} \text{ for } k=1, \dots, k_0-1$, we refer to $\Delta^{(1)}_{\mathbf{p}}$ as a \textit{standard probabilistic Laplacian}. A direct computation shows that the spectral decimation function corresponding to a standard probabilistic Laplacian can be derived from a Chebyshev polynomial. Recall that the Chebyshev polynomials of the first kind satisfy the   recurrence relations
%\begin{align}
%\label{eq:recurrenceRelationsCheb}\begin{cases}
 $   Cheb_0(z)=1$, $Cheb_1(z) = z$, $    Cheb_k(z) = 2z$, $ Cheb_{k-1}(z) - Cheb_{k-2}(z) $, where $ k \in \{2,\dots, k_0\}$.
%  \end{cases}
%\end{align}
%
%\begin{proposition}
%\label{prop:SpecDeciChebeshev}
If $\Delta^{(1)}_{\mathbf{p}}$ is a standard probabilistic Laplacian, i.e. given by (\ref{eq:initialLapNandD})
and  $ p_k=\frac{1}{2} \text{ for } k=1, \dots, k_0-1$, then the corresponding spectral decimation function in (\ref{eq:SpecDeciFuncForSeqOfLap})  is $R_{\Delta^{(1)}_{\mathbf{p}}}(z)=1-  (-1)^{k_0} Cheb_{k_0}(z-1)$.
%\begin{align*}
%R_{\Delta^{(1)}_{\mathbf{p}}}(z)=1-  (-1)^{k_0} Cheb_{k_0}(z-1).
%\end{align*}
%\end{proposition}
This is discussed in detail in \cite[Section 4]{BajorinVibrationSpectra2008}. 
It is evident that for $k_0=2$ the $3\times3$ standard probabilistic Laplacian  is the only possible centrosymmetric probabilistic Laplacian of this size. 

When $k_0=3$,  then $\Delta^{(1)}_{\mathbf{p}}$ is of size $4\times4$ with $\mathbf{p}=\{1,p_1, p_2, 0\}$. Imposing the centrosymmetry assumption on the transition probabilities gives $p_2=1-p_1$. Hence, we obtain a one-parameter family of centrosymmetric probabilistic Laplacians and it is not difficult to see that  this type of probabilistic Laplacians was already studied in \cite{TeplyaevSpectralZeta2007} under the terminology of a \textit{pq-model}, see also \cite{Ulysses2017,ChenTeplyaev2016}   and more recent \cite[Proposition 2.3]{BaluMogOkoTep2021spectralAMO}.

\subsection{One-parameter family of   Laplacians  with $k_0=4$}

Now we consider the case $k_0=4$, i.e. a centrosymmetric probabilistic Laplacian $\Delta^{(1)}_{\mathbf{p}}$ of size $5\times5$ with $\mathbf{p}=\{1,p_1, p_2,p_3, 0\}$. Imposing the centrosymmetry assumption on the transition probabilities gives $p_3=1-p_1$ and $p_2=1-p_2$.  We obtain again a one-parameter family of centrosymmetric probabilistic Laplacians whose matrices are given by
 
\begin{align}
\label{eq:5x5mirrorSymProbLap}
\Delta^{(1)}_{\mathbf{p}} =
\left(\begin{matrix}1 & -1 & 0 & 0 & 0\\p - 1 & 1 & - p & 0 & 0\\0 & - \frac{1}{2} & 1 & - \frac{1}{2} & 0\\0 & 0 & - p & 1 & p - 1\\0 & 0 & 0 & -1 & 1\end{matrix}\right), 
\quad
\Delta^{(1),D}_{\mathbf{p}} =
\left(\begin{matrix}
 1 & - p & 0 \\
 - \frac{1}{2} & 1 & - \frac{1}{2} \\
 0 & - p & 1 
\end{matrix}\right)
\end{align}
where we set $p_1=p$ for some $p \in (0,1)$. Note that this type of probabilistic Laplacians is related to the graph Laplacians studied in \cite{FangKingStrichartz2019}.
\begin{figure}
  \begin{minipage}[b]{1.\linewidth}
    \centering
    \includegraphics[width=0.75\linewidth]{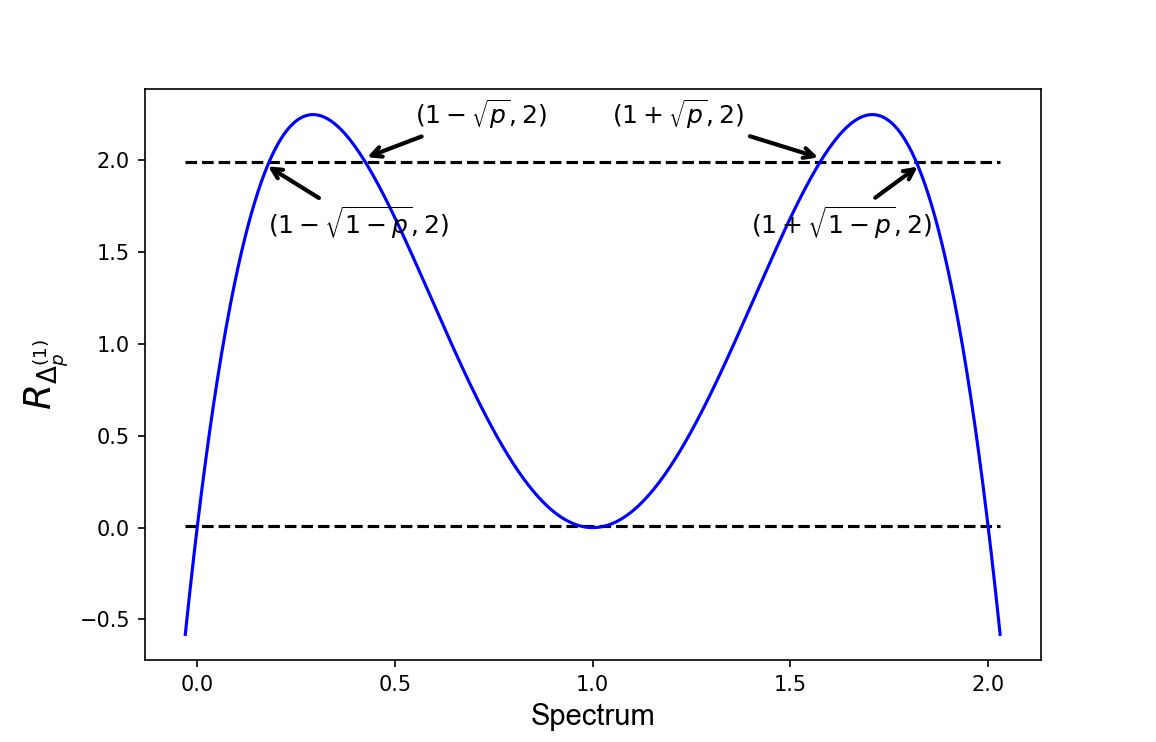} 
    \end{minipage}%%
    \\
      \begin{minipage}[b]{1.\linewidth}
    \centering
\includegraphics[width=0.7\linewidth]{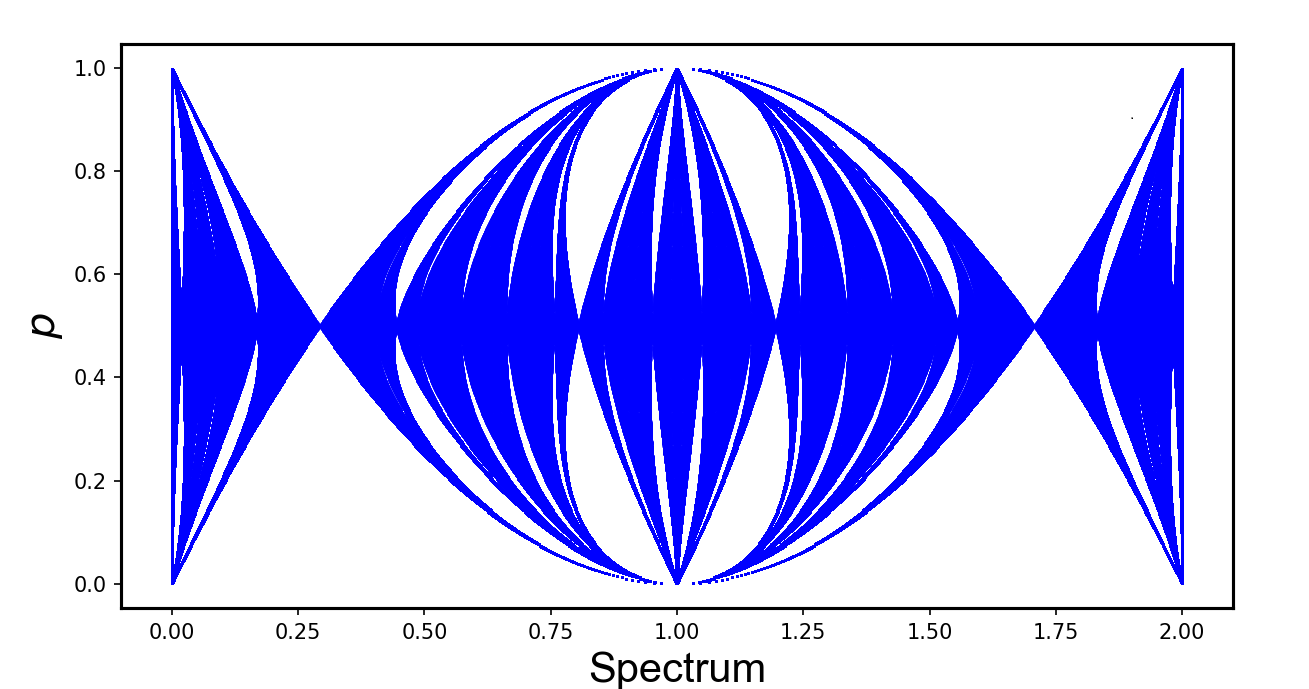} 
  \end{minipage} 
  \caption{ (Top) The spectral decimation function (\ref{eq:SpecDecik0eq4}) corresponding to the probabilistic Laplacians (\ref{eq:5x5mirrorSymProbLap}). The preimages in (\ref{eq:preImages}) are indicated using the cutoffs $y=0$ and $y=2$. (Bottom) The spectrum of $\Delta^{(\ell)}_{\mathbf{p}}$ for $\ell=6$ (x-axis) is plotted for $p \in (0, 1)$ (y-axis). }
  \label{fig:pqModelBobStud}
\end{figure}
Direct computation gives
\begin{align}
%\label{eq:Spectrak0eq3NandD}
\begin{cases}
\sigma \big( \Delta^{(1)}_{\mathbf{p}} \big)=\Big\{ 0 , \  1 - \sqrt{1 - p}, \  1, \  1+\sqrt{1 - p} , \  2\Big\} ,\\ \sigma \big( \Delta^{(1),D}_{\mathbf{p}} \big) = \Big\{  1 - \sqrt{p} , \ 1 ,\  1+\sqrt{p} \Big\}.
\end{cases}
\end{align}
We compute the spectral decimation function using the formula in  \eqref{eq:SpecDeciFuncForSeqOfLap},
\begin{align}
\label{eq:SpecDecik0eq4}
R_{\Delta^{(1)}_{\mathbf{p}}}(z) =   \frac{2 z \left(2- z \right) \left(z - 1\right)^{2}}{p \left( 1-p\right)}.
\end{align}
We obtain for the set of the preimages
\begin{align}
\label{eq:preImages}
\begin{cases}
  R^{-1}_{\Delta^{(1)}_{\mathbf{p}}}(0)=\Big\{0,1,2\Big\}, \\ R^{-1}_{\Delta^{(1)}_{\mathbf{p}}}(2)=\Big\{    1 - \sqrt{1 - p}, \ 1 - \sqrt{p} ,\  1+ \sqrt{p} , \  1+ \sqrt{1 - p}  \Big\}.
  \end{cases}
\end{align}
Figure \ref{fig:pqModelBobStud} (top) shows the preimages in (\ref{eq:preImages})  using the cutoffs $y=0$ and $y=2$. The sequence of self-similar probabilistic Laplacians $\Delta^{(\ell)}_{\mathbf{p}}$ is then constructed as in Definition \ref{def:selfSimiLapSeq}, where $\Delta^{(1)}_{\mathbf{p}}$ is initialized by \eqref{eq:5x5mirrorSymProbLap}. Figure \ref{fig:pqModelBobStud} (bottom) shows how the spectrum of $\Delta^{(\ell)}_{\mathbf{p}}$, $\ell=6$ changes when varying the parameter $p \in (0,1)$.

\subsection{Two-parameter family of   Laplacians with $k_0=5$}

Now we consider the case $k_0=5$, i.e. a centrosymmetric probabilistic Laplacian $\Delta^{(\ell)}_{\mathbf{p}}$ of size $6\times6$ with $\mathbf{p}=\{1,p_1, p_2,p_3, p_4,0\}$. Imposing the centrosymmetry assumption on the transition probabilities gives $p_4=1-p_1$ and $p_3=1-p_2$. We obtain in this case a two-parameter family of centrosymmetric probabilistic Laplacians whose matrices are given by
{\fontsize{8}{10} 
\begin{align}
\label{eq:2ParameterLap}
\Delta^{(1)}_{\mathbf{p}} =
\left(\begin{matrix}1 & -1 & 0 & 0 & 0 & 0\\
p_{1} - 1 & 1 & - p_{1} & 0 & 0 & 0\\
0 & p_{2} - 1 & 1 & - p_{2} & 0 & 0\\
0 & 0 & - p_{2} & 1 & p_{2} - 1 & 0\\
0 & 0 & 0 & - p_{1} & 1 & p_{1} - 1\\
0 & 0 & 0 & 0 & -1 & 1\end{matrix}\right),
\
\Delta^{(1),D}_{\mathbf{p}} =
\left(\begin{matrix}
 1 & - p_{1} & 0 & 0 \\
 p_{2} - 1 & 1 & - p_{2} & 0 \\
 0 & - p_{2} & 1 & p_{2} - 1 \\
 0 & 0 & - p_{1} & 1 \end{matrix}\right).
\end{align} }%
We compute the spectral decimation function using the formula in  (\ref{eq:SpecDeciFuncForSeqOfLap}),{\fontsize{8}{10} 
\begin{align*}
R_{\Delta^{(1)}_{\mathbf{p}}}(z) =
\frac{z \left(z^{2} +z (p_{2}   - 3) + p_{1} p_{2} - 2 p_{2}  + 2\right) \left(  z^{2} - z(p_{2} +2)  +  p_{1} p_{2} - p_{1} + p_{2}   + 1\right)}{p_{1} p_{2} \left(p_{1} - 1\right) \left(p_{2} - 1\right)}
\end{align*}}%

\begin{figure}[htb]
\centering
%\hspace*{-2cm}
\includegraphics[width=0.7\textwidth]{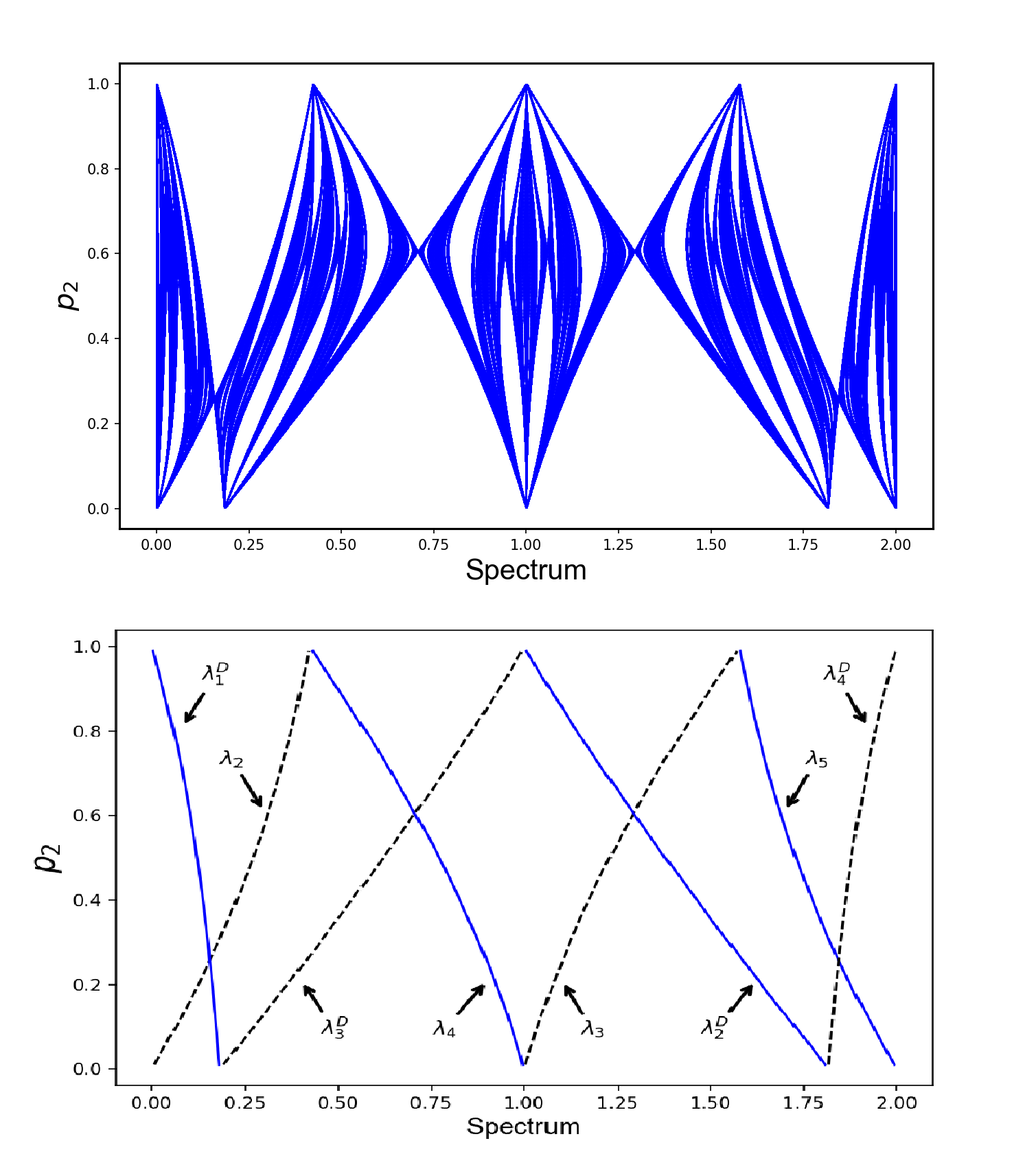}
%\vspace{0cm}
  \caption{(Top) The spectrum of $\Delta^{(\ell)}_{\mathbf{p}}$ for $\ell=5$ is plotted on the x-axis for $p_2 \in (0, 1)$ (y-axis), while $p_1=\frac{2}{3}$ is fixed. (Bottom) Plotting the eigenvalues as functions of $p_2$ reproduces the contour of the spectra.}
  \label{spiderp1eq2over3Varyp2WithBoundaries500Level5}
\end{figure}

For the convenience of the notation, we fix $p_1 = \frac{2}{3}$ and compute the eigenvalues as functions of the parameter $p_2$,{\fontsize{9}{10} 
\begin{align}
\label{eq:eigenvaluesk0eq5}
\sigma \big( \Delta^{(1)}_{\mathbf{p}} \big)=\Big\{ \lambda_1 , \  \lambda_2, \  \lambda_3, \  \lambda_4 , \  \lambda_5, \ \lambda_6\Big\} ,\quad \sigma \big( \Delta^{(1),D}_{\mathbf{p}} \big) = \Big\{  \lambda^D_1 , \  \lambda^D_2 ,\  \lambda^D_3, \  \lambda^D_4 \Big\}.
\end{align}}%
The formulas for eigenvalues are given by,
{\fontsize{8}{10} 
\begin{align*}
 \lambda_1 =0,  \quad  \lambda_2 =\frac{1+p_{2}}{2} - \frac{ \sqrt{9 p_{2}^{2} - 6 p_{2} + 9}}{6} ,      \quad  \lambda_3 =\frac{1+p_{2}}{2} + \frac{ \sqrt{9 p_{2}^{2} - 6 p_{2} + 9}}{6} ,                  \\
\lambda_4 = \frac{3-p_{2}}{2} - \frac{ \sqrt{9 p_{2}^{2} - 6 p_{2} + 9}}{6}, \quad  \lambda_5 = \frac{3-p_{2}}{2} + \frac{ \sqrt{9 p_{2}^{2} - 6 p_{2} + 9}}{6} , \quad  \lambda_6=2,
\end{align*}
and
\begin{align*}
 \lambda^D_1= \frac{2-p_{2}}{2} - \frac{ \sqrt{9 p_{2}^{2} - 24 p_{2} + 24}}{6}  , \quad \lambda^D_2= \frac{2-p_{2}}{2} + \frac{ \sqrt{9 p_{2}^{2} - 24 p_{2} + 24}}{6}, \\
 \lambda^D_3=\frac{2+p_{2}}{2} - \frac{ \sqrt{9 p_{2}^{2} - 24 p_{2} + 24}}{6}, \quad \lambda^D_4=\frac{2+p_{2}}{2} + \frac{ \sqrt{9 p_{2}^{2} -24 p_{2} + 24}}{6}.
\end{align*} }%

Now we construct the sequence of self-similar probabilistic Laplacians $\Delta^{(\ell)}_{\mathbf{p}}$ as in Definition \ref{def:selfSimiLapSeq}, where $\Delta^{(1)}_{\mathbf{p}}$ is initialized by \eqref{eq:2ParameterLap}. 
Figure \ref{spiderp1eq2over3Varyp2WithBoundaries500Level5} (top) shows how the spectrum of $\Delta^{(\ell)}_{\mathbf{p}}$, $\ell=5$ changes when varying the parameter $p_2 \in (0,1)$ and fixing $p_1 = \frac{2}{3}$. Plotting the eigenvalues (\ref{eq:eigenvaluesk0eq5}) as functions of $p_2$ reproduces the contour of the spectra in Figure \ref{spiderp1eq2over3Varyp2WithBoundaries500Level5}(top). Using the eigenvalue-formulas, we can for example determine when the spectral gap between $\lambda_1^D$ and $\lambda_2$ closes, namely by reducing the problem to solving the equation $\lambda_1^D=\lambda_2$, see Figure \ref{spiderp1eq2over3Varyp2WithBoundaries500Level5} (bottom).

\subsection*{Acknowledgments}

The work of G. Mograby and K. Okoudjou was supported by ARO grant W911NF1910366.  K. Okoudjou was additionally supported by NSF DMS-1814253. A.~Teplyaev was partially supported by NSF DMS grant 1613025 and by the Simons Foundation.

\bibliographystyle{plain}
\bibliography{BibList-lastUpdate-4Aug2021}

\end{document}